\documentclass[10pt,a4paper,english]{article}
\usepackage[a4paper,left=3cm,right=2cm]{geometry}
\usepackage{babel}
\usepackage[latin1]{inputenc}
\usepackage{amsmath}
\usepackage{amsthm}
\usepackage{amsfonts}
\usepackage{indentfirst}
\usepackage{graphicx}
\newtheorem{de}{Definition}
\newtheorem{pro}{Proposition}
\newtheorem{cor}{Corollary}
\newtheorem{teo}{Theorem}
\newtheorem{rem}{Remark}
\newtheorem{lem}{Lemma}
\newtheorem{exa}{Example}
\newtheorem{alg}{Algorithm}

\newcommand{\co}{{\mathcal O}}
\newcommand{\ci}{{\mathcal I}}

\newcommand{\gp}{\mathbb{P}}

\newcommand{\gz}{\mathbb{Z}}
\newcommand{\gq}{\mathbb{Q}}
\newcommand{\gr}{\mathbb{R}}
\newcommand{\gc}{\mathbb{C}}

\newcommand{\ck}{\mathcal K}
\newcommand{\cb}{{\mathcal B}}

\newcommand{\cp}{{\mathcal P}}

\newcommand{\cf}{{\mathcal F}}
\newcommand{\cd}{{\mathcal D}}

\newcommand{\findemo}{$\ \ \square$}

\newcommand{\card}{{\rm card}}
\newcommand{\con}{{\rm con}}

\renewcommand{\int}{{\rm int}}

\newcommand{\pic}{{\rm Pic}}

\usepackage{color}

\title{ {\bf Algebraic Integrability of Foliations of the Plane}}
\author{C. Galindo \thanks{Supported by the Spain Ministry of Education
 MTM2004-00958, GV05/029  and Bancaixa P1-1A2005-08
}
 \and F.  Monserrat$^{*}$}
\date{}

\begin{document}
\maketitle

\begin{abstract}
We give an algorithm to decide whether an algebraic plane
foliation $\cf$ has a rational first integral and to compute it in
the affirmative case. The algorithm runs whenever we assume the
polyhedrality of the cone of curves of the surface  obtained after
blowing-up the set ${\cal B}_{\cal F}$ of infinitely near points
needed to get the dicritical exceptional divisors of a minimal
resolution of the singularities of $\cf$. This condition can be
detected in several ways, one of them from the proximity relations
in ${\cal B}_{\cal F}$ and, as a particular case, it holds when
the cardinality of ${\cal B}_{\cal F}$ is less than 9.
\end{abstract}
\section{Introduction}

The problem of deciding whether  a complex polynomial differential
equation on the plane is algebraically integrable goes back to the
end of the nineteenth century when  Darboux \cite{dar}, Poincar\'e
\cite{poi1, poi2, poi3},  Painlev\'e \cite{pai} and Autonne
\cite{aut} studied it. In modern terminology and from a more
algebraic point of view, it can be stated as deciding whether an
algebraic foliation $\cal{F}$ with singularities on the projective
plane over an algebraically closed field of characteristic zero
(plane foliation or foliation on $\gp^2$, in the sequel) admits a
rational first integral and, if it is so, to compute it. In this
paper, we shall give a satisfactory answer to that problem when
the cone of curves of certain surface  is polyhedral. This surface
is obtained by blowing-up what we call dicritical points of a
minimal resolution of the singularities of $\cf$.

The existence of a rational first integral is  equivalent to state
that every invariant curve of $\cal{F}$ is algebraic. The fact
that $\cal{F}$ admits algebraic invariant curves has interest for
several reasons. For instance, it is connected with the center
problem for quadratic vector fields \cite{sch, ch-gi}, with
problems related to solutions of Einstein's field equations  in
general relativity \cite{hew} or with the second part of the
Hilbert's sixteenth problem \cite{hil} (see also \cite{sma}),
which looks for a bound of the number of limit cycles for a (real)
polynomial  vector field
(see for example \cite{lli} and precedents in the proofs of
results in \cite{bau}).

Coming back to the nineteenth century and in the analytic complex
case, it was Poincar\'e \cite{poi2} who observed that ``to find out
if a differential  equation of the first order and of the first
degree is algebraically integrable, it is enough to find an upper
bound for the degree of the integral. Afterwards, we  only need to
perform purely algebraic computations". This observation gave rise
to the so called Poincaré problem which, nowadays, is established as
the one of bounding the degrees of the algebraic leaves  of a
foliation whether  it is algebraically integrable or not. It was
Poincar\'e himself who studied a particular case within the one
where the singularities of the foliation are non-degenerated
\cite{poi2}.
Carnicer \cite{car} provided an answer for the nondicritical
foliations. Bounds depending on the invariant curve have been
obtained by Campillo and Carnicer \cite{ca-ca} and, afterwards,
improved by Esteves and Kleiman \cite{es-kl}. However, the classical
Poincar\'e problem has a negative answer, that is the degree of a
general irreducible invariant curve of a plane foliation $\cal{F}$
with a rational first integral cannot be bounded by a function on
the degree of the foliation, and this happens even in the case of
families of foliations where the analytic type of each singularity
is constant \cite{l-n}. An analogous answer is given in \cite{l-n}
for a close question posed by Painlev\'e in \cite{pai}, which
consists of recognizing the genus of the general solution of a
foliation as above.

On the other hand, Prelle and Singer gave in \cite{pr-si} a
procedure to compute elementary first integrals of foliations
$\cal{F}$ on the projective plane over the complex numbers. As a
particular case, it uses results by Darboux and Jouanolou to deal
with the computation of meromorphic first integrals of $\cal{F}$;
however, the obstruction revealed by the Poincar\'e problem makes
the above procedure be only a semi-decision one (see the
implementation by Man in \cite{man1} and that given in \cite{man2}
for the rational case). Notice that Man uses packages involving
Groebner basis to detect inconsistency as well as to solve
consistent systems of equations \cite[3.3]{man1}. In this paper, we
shall give an alternative algorithm (Algorithm \ref{alg1}), uniquely
involving the resolution of   systems of linear equations that, for
each tentative degree $d$ of a general irreducible invariant curve,
decides whether $\cf$ has a rational first integral of degree $d$
and computes it in the affirmative case.

The so called $d$-extactic curves of a plane foliation $\cal{F}$,
studied by Lagutinskii \cite{do-lo} and Pereira \cite{per},
provide a nice, but non-efficient, procedure to decide whether
$\cal{F}$ has a rational first integral of some given degree.
Moreover, two sufficient conditions so that $\cf$ had a rational
first integral are showed in \cite{cha-lli}.


The main result of this paper is to give an algorithm to decide
whether  plane foliations in a certain class have a rational first
integral, which also allows to compute it in an effective manner.
In fact, in the affirmative case, we obtain a primitive first
integral, that allows to get any other first integral, and,
trivially, a bound for the degree of the irreducible components of
the invariant curves of the foliation. A well-known result is the
so called resolution theorem \cite{be, le, seid} that asserts that
after finitely many blow-ups at singular points (of the
successively obtained foliations by those blow-ups), $\cal{F}$ is
transformed in a foliation on another surface with finitely many
singularities, all of them of a non-reducible by blowing-up type,
called simple. The input of our algorithm will be the foliation
and also a part of the mentioned configuration of infinitely near
points that resolves its singularities.
This part is what leads to obtain the so called dicritical
exceptional divisors of the foliation. If we blow  that
configuration, we get a surface $Z_\cf$, and the cone of curves
$NE(Z_\cf)$ of that surface will be our main tool (such a cone is a
basic object in the minimal model theory \cite{kollar}). When $\cf$
admits a rational first integral, it has associated an irreducible
pencil of plane curves whose general fibers provide a divisor
$D_\cf$ on $Z_\cf$. The class in the Picard group of $D_\cf$
determines a face of $NE(Z_\cf)$ and it has codimension 1 if, and
only if, $\cal{F}$ has an independent system of algebraic solutions
$S$ (see Definition \ref{condition}). If we get a system as $S$,
then we can determine   $D_\cf$ (see Theorem \ref{gordo2} and
Proposition \ref{memo}). We devote Section \ref{DOS} to explain and
prove the above considerations and Section \ref{UNO} to give the
preliminaries and notations.

In  Section \ref{TRES}, we show our main result by means of two
algorithms, based in the above results, that one must jointly use:
Algorithm \ref{alg3} runs for foliations $\cal{F}$ whose cone of
curves $NE(Z_\cf)$ is (finite) polyhedral and computes a system $S$
as above (or discards the existence of a rational first integral),
whereas Algorithm \ref{alg2} uses $S$ to compute a rational first
integral (or newly discards its existence). Both algorithms can be
implemented without difficulty. Polyhedrality of the cone of curves
happens in several cases. For instance, when the anti-canonical
bundle on $Z_\cf$ is ample by the Mori Cone Theorem \cite{mor}. Also
whenever $Z_\cf$ is obtained by blowing-up configurations of toric
type \cite{cam, oda}, relative to pencils $H - \lambda Z^d$, where
$H$ is an homogeneous polynomial of degree $d$ that defines a curve
with a unique branch at infinity \cite{camp}, or of less than nine
points \cite{cam-gon} (see also \cite{man})
that are included within the much wider set of P-sufficient
configurations (Definition \ref{suficiente}), whose cone of curves
is also polyhedral (see \cite{gal1} for a proof). Notice that in
this last case (P-sufficient configurations) and when Algorithm
\ref{alg3} returns an independent system of algebraic solutions, to
compute a rational first integral of $\cf$ is simpler than above
and, furthermore, the Painlev\'e problem is solved (Propositions
\ref{memo} and \ref{lacinco}). Finally, this last section includes
several illustrative examples of  our results.


\section{Preliminaries and notations}\label{UNO}


Let $k$ be an algebraically closed field  of characteristic zero. An
(algebraic singular) {\it foliation} $\cf$ on a projective smooth
surface (a surface in the sequel) $X$ can be defined by the data
$\{(U_i,\omega_i)\}_{i\in I}$, where $\{U_i\}_{i\in I}$ is an open
covering of $X$, $\omega_i$ is a non-zero regular differential
1-form on $U_i$ with isolated zeros and, for each couple $(i,j)\in
I\times I,$
\begin{equation}\label{eq}
\omega_i=g_{ij}w_j \mbox{\;\;\;on \;\;} U_i \cap U_j,
\mbox{\;\;\;\;} g_{ij}\in \co_X(U_i\cap U_j)^*.
\end{equation}


Given $p\in X$, a {\it (formal) solution of $\cf$ at $p$} will be
an irreducible element $f\in \widehat{\co}_{X,p}$ (where
$\widehat{\co}_{X,p}$ is the ${\rm m}_p$-adic completion of the
local ring $\co_{X,p}$ and ${\rm m}_p$ its maximal ideal) such
that the local differential 2-form $\omega_p \wedge df$ is a
multiple of $f$, $w_p$ being a local equation of $\cf$ at $p$. An
element in $\widehat{\co}_{X,p}$ will be said to be {\it invariant
by} $\cf$ if all its irreducible components are solutions of $\cf$
at $p$. An {\it algebraic solution} of $\cf$ will be an integral
(i.e., reduced and irreducible) curve $C$ on $X$ such that its
local equation at each point in its support is invariant by $\cf$.
Moreover, if every integral component of a curve $D$ on $X$ is an
algebraic solution, we shall say that $D$ {\it is invariant by}
$\cf$.

The transition functions $g_{ij}$ of a foliation $\cf$ define an
invertible sheaf $\cal N$ on $X$, the {\it normal sheaf} of $\cf$,
and the relations (\ref{eq}) can be thought as defining relations
of a global section of the sheaf ${\cal N}\otimes \Omega_X^1$,
which has isolated zeros (because each $\omega_i$ has isolated
zeros). This section is uniquely determined by the foliation
$\cf$, up to multiplication by a non zero element in $k$.
Conversely, given an invertible sheaf ${\cal N}$ on $X$, any
global section of ${\cal N}\otimes \Omega_X^1$ with isolated zeros
defines a foliation $\cf$ whose normal sheaf is ${\cal N}$. So,
alternatively, we can define a foliation as a map of
$\mathcal{O}_{X}$-modules $\cf:\Omega^{1}_X \rightarrow {\cal N}$,
${\cal N}$ being some invertible sheaf as above. Set $Sing(\cf)$
the  singular locus of $\cf$, that is the subscheme of $X$ where
$\cf$ fails to be surjective.

Particularizing to the projective plane $\mathbb{P}_k^2 = \gp^2$,
for a non-negative integer $r$, the Euler sequence, $ 0 \rightarrow
\Omega^1_{\gp^2} \rightarrow \mathcal{O}_{\gp^2}(-1)^3 \rightarrow
\mathcal{O}_{\gp^2} \rightarrow 0 $, allows to regard the foliation
$\cf: \Omega^1_{\gp^2} \rightarrow \mathcal{O}_{\gp^2}(r-1)$, in
analytic terms, as induced by a homogeneous vector field $\mathbf{X}
= U \partial / \partial X + V
\partial / \partial Y + W \partial / \partial Z$, where $U,V,W$ are
homogeneous polynomials of degree $r$ in homogeneous coordinates
$(X:Y:Z)$ on $\gp^2$. It is convenient to notice that two vector
fields define the same foliation if, and only if, they differ by a
multiple of the radial vector field. By convention, we shall say
that $\cf$ has degree $r$. We prefer to use forms to treat
foliations and, so, equivalently, we shall give a foliation $\cf$
on $\gp^2$ of degree $r$,
up to a scalar factor,
by means
of a projective 1-form
$$
\Omega=A dX + B dY + C dZ,
$$
where $A,B$ and $C$ are homogeneous polynomials of degree $r+1$
without common factors which satisfy the Euler's condition
$XA+YB+ZC=0$ (see \cite{G-M}). $\Omega$ allows to handle easily the
foliation in local terms and the singular points of $\cf$ are the
common zeros of the polynomials $A, B$ and $C$.
Moreover, a curve $D$ on
$\gp^2$ is invariant by $\cf$ if, and only if, $G$ divides the
projective 2-form $dG \wedge \Omega$, where $G(X:Y:Z)=0$ is an
homogeneous equation of $D$.

To blow-up a surface at a closed point and the corresponding
evolution of a foliation on it, will be an important tool in this
paper. Thus, let us consider a sequence of morphisms
\begin{equation}
\label{seq} X_{n+1} \mathop  {\longrightarrow} \limits^{\pi _{n} }
X_{n} \mathop {\longrightarrow} \limits^{\pi _{n-1} }  \cdots
\mathop {\longrightarrow} \limits^{\pi _2 } X_2 \mathop
{\longrightarrow} \limits^{\pi _1 } X_1 : = \gp^2,
\end{equation}
where $\pi_i$ is the blow-up of $X_i$ at a closed point $p_i\in
X_i$, $1\leq i\leq n$. The associated set of closed points $ \ck =
\{p_1,p_2,\ldots,p_n\}$ will be called a {\it configuration} over
$\gp^2$ and the variety $X_{n+1}$ the {\it sky} of $\ck$; we
identify two configurations with $\gp^2$-isomorphic skies. We
shall denote by $E_{p_i}$ (respectively, $\tilde{E}_{p_i}$)
(respectively, $E_{p_i}^*$) the exceptional divisor appearing in
the blow-up $\pi_i$ (respectively, its strict transform on
$X_{n+1}$) (respectively, its total transform on $X_{n+1}$).
Also, given two points $p_i, p_j$ in $\ck$, we shall say that $p_i$
is {\it infinitely near to} $p_j$ (denoted $p_i \geq p_j$) if either
$p_i=p_j$ or $i>j$ and $\pi_j\circ \pi_{j+1}\circ \cdots \circ
\pi_{i-1}(p_i)=p_j$. The relation $\geq$ is a partial ordering among
the points of the configuration $\ck$. Furthermore, we say that a
point $p_i$ is {\it proximate} to other one $p_j$ whenever $p_i$ is
in the strict transform of the exceptional divisor created after
blowing up at $p_j$ in the surface which contains $p_i$. As a visual
display of a configuration $\ck$, we shall use the so called {\it
proximity graph} of $\ck$,  whose vertices represent those points in
$\ck$, and two vertices, $p,q \in \ck$, are joined by an edge if $p$
is proximate to $q$. This edge is dotted except when $p$ is in the
first infinitesimal neighborhood of $q$ (here the edge is
continuous). For simplicity sake, we delete those edges which can be
deduced from others.

If $\cf$ is a foliation on $\gp^2$, the sequence of morphisms
(\ref{seq}) induces, for each $i=2,3,\ldots,n+1$, a foliation
$\cf_i$ on $X_i$, the strict transform of $\cf$ on $X_i$ (see
\cite{brun}, for instance).
As we have said in the Introduction, Seidenberg's result of
reduction of singularities \cite{seid} proves  that there is a
sequence of blow-ups as in (\ref{seq}) such that the strict
transform $\cf_{n+1}$ of $\cf$ on the last obtained surface
$X_{n+1}$ has only certain type of singularities which cannot be
removed by blowing-up, called simple singularities. Such a
sequence of blow-ups is called a {\it resolution} of $\cf$, and it
will be {\it minimal} if it is so with respect to the number of
involved blow-ups. Assuming that (\ref{seq}) is a minimal
resolution of $\cf$, we shall denote by $\ck_{\cf}$ the associated
configuration $\{p_i\}_{i=1}^n$. Note that each point $p_i$ is an
ordinary (that is, not simple) singularity  of the foliation
$\cf_i$.

Dicriticalness of divisors and points will be an essential concept
to decide if certain plane foliations have a rational first
integral, object of our study. Next we state the definitions.

\begin{de} {\rm An exceptional divisor $E_{p_i}$ (respectively, a point $p_i\in
\ck_{\cf}$) of a minimal resolution of a plane foliation $\cf$ is
called {\it non-dicritical} if it is invariant by the foliation
$\cf_{i+1}$ (respectively, all the exceptional divisors $E_{p_j}$,
with $p_j\geq p_i$, are non-dicritical). Otherwise, $E_{p_i}$
(respectively,  $p_i$) is said to be {\it dicritical}.}
\end{de}

Along this paper, we shall denote by ${\cal B}_{\cal F}$
the configuration of dicritical points in ${\cal K}_{\cal F}$, and
by ${\cal N}_{\cal F}$  the set of points $p_i\in {\cal B}_{\cal F}$
such that $E_{p_i}$ is a non-dicritical exceptional divisor.


\begin{de}
{\rm We shall say that a plane foliation $\cf$ has a {\it rational
first integral} if there exists a rational function $R$ of $\gp^2$
such that $dR \wedge \Omega=0$. }
\end{de}

The existence of a rational first integral is equivalent to each one
of the following three facts (see \cite{jou}): $\cf$ has infinitely
many  algebraic solutions, all the solutions of $\cf$ are
restrictions of algebraic solutions and there exists a unique
irreducible pencil of plane curves $\cp_{\cf}:=\left\langle {F,G}
\right\rangle\subseteq H^0(\gp^2,\co_{\gp^2}(d))$, for some $d\geq
1$, such that the algebraic solutions of $\cf$ are exactly the
integral components of the curves of the pencil. {\it Irreducible
pencil} means that its general elements are integral curves.



Two generators $F$ and $G$ of $\cp_{\cf}$ give rise to a rational
first integral $R=\frac{F}{G}$ of $\cf$ and, if $T$ is whichever
rational function of $\gp^1$, then $T(R)$ is also a rational first
integral of $\cf$ and any rational first integral is obtained in
this way.
We shall consider rational first integrals arising from the
unique pencil $\cp_{\cf}$.

Until the end of this section, we shall assume that $\cf$ has a
rational first integral. $\cp_{\cf}$ has finitely many base points,
since $F$ and $G$ have no common factor.
Set $\ci\subseteq \co_{\gp^2}$  the ideal sheaf supported at the
base points of $\cp_{\cf}$ and such that $\ci_p=(F_p, G_p)$ for each
such a point $p$, where $F_p$ and $G_p$ are the natural images of
$F$ and $G$ in $\co_{\gp^2,p}$. There exists a sequence of blow-ups
centered at closed points
\begin{equation}
\label{seq2} X_{m+1} \mathop  {\longrightarrow} \limits^{\pi _{m}
} X_{m} \mathop {\longrightarrow} \limits^{\pi _{m-1} }  \cdots
\mathop {\longrightarrow} \limits^{\pi _2 } X_2 \mathop
{\longrightarrow} \limits^{\pi _1 } X_1 : = \gp^2
\end{equation}
such that, if $\pi_{\cf}$ denotes the composition morphism
$\pi_1\circ \pi_2\circ \cdots \pi_m$, the pull-back $\pi_{\cf}^*
\ci$ becomes an invertible sheaf of $X_{m+1}$ \cite{cam}. We
denote by ${\cal C}_{\cf}$ the set of  centers of the blow-ups
that appear in a minimal sequence with this property and by
$Z_\cf$ the sky of ${\cal C}_{\cf}$. This sequence can also be
seen as a minimal sequence of blow-ups that eliminate the
indeterminacies of the rational map $\gp^2 \cdots \rightarrow
\gp^1$ induced by the pencil $\cp_{\cf}$. Hence, there exists a
morphism $h_{\cf}:Z_{\cf}\rightarrow \gp^1$ factorizing through
$\pi_{\cf}$. Notice that this morphism is essentially unique, up
to composition with an automorphism of $\gp^1$.

Let $F$ and $G$ be two general elements of the pencil $\cp_{\cf}$
and, for each $p_i\in \ck_{\cf}$, assume that
$f$ (respectively, $g$) gives a local equation at $p_i$ of the
strict transform of the curve on $\gp^2$ defined by $F$
(respectively, $G$). Then, the local solutions of $\cf_i$ at $p_i$
are exactly the irreducible components of the elements of the
(local) pencil in $\widehat{\co}_{X_i,p_i}$ generated by $f$ and
$g$ \cite{julio}. As a consequence, the following result is clear


\begin{pro}
If ${\cal F}$ is a foliation on $\gp^2$ with a rational first
integral, then the configurations ${\cal C}_{\cal F}$ and ${\cal
B}_{\cal F}$ coincide.

\end{pro}

\section{Foliations with a rational first integral} \label{DOS}

Along this paper, we shall consider a foliation ${\cal F}$ on
$\gp^2$ such that its associated configuration $\cb_{\cf}$ has
cardinality larger than 1 and we keep  the notations as in the
previous section.
Notice that in case that ${\cal B}_{\cf}$ be empty, it is obvious
that the foliation $\cf$ has no rational first integral and, when
${\cal B}_{\cf}$ consists of a point, only quotients of linear
homogeneous polynomials defining transversal lines which pass
through that point can be rational first integrals.


From now on,   $Z_{\cf}$ will denote the sky of $\cb_{\cf}$ (it will
be the one of ${\cal C}_{\cf}$ whenever $\cf$ has a rational first
integral). Denote by $A(Z_{\cf})$ the real vector space (endowed
with the usual real topology) $\pic(Z_{\cf})\otimes_{\gz} \gr\cong
\gr^{m+1}$, where $\pic(Z_{\cf})$ stands for the Picard group of the
surface $Z_{\cf}$. The {\it cone of curves} (respectively, {\it nef
cone}) of $Z_{\ck}$, which we shall denote by $NE(Z_{\cf})$
(respectively, $P(Z_{\cf})$), is defined to be  the convex cone of
$A(Z_{\cf})$
 generated by the images of the effective (respectively, nef) classes
 in $\pic(Z_{\cf})$. The $\gz$-bilinear form  $\pic(Z_{\cf})\times
\pic(Z_{\cf})\rightarrow \gz$ given by Intersection Theory induces a
non-degenerate $\gr$-bilinear pairing
\begin{equation}
\label{pairing} A(Z_{\cf}) \times A(Z_{\cf}) \rightarrow \gr.
\end{equation}
For each pair  $(x,y)\in A(Z_{\cf})\times A(Z_{\cf})$, $x\cdot y$
will denote its image by the above bilinear form.

On the other hand, given  a convex cone $C$ of $A(Z_{\cf})$, its
{\it dual cone} is defined to be $C^\vee:=\{x\in A(Z_{\cf}) \mid
x\cdot y\geq 0 \;\;\mbox{for all}\;\;y\in C\}$, and a {\it face}
of $C$ is a sub-cone $D\subseteq C$ such that $a+b\in D$ implies
that $a,b\in D$, for all pair of elements $a,b\in C$. The
$1$-dimensional faces of $C$ are also called {\it extremal rays}
of $C$. Note that $P(Z_{\cf})$ is the dual cone of $NE(Z_{\cf})$,
and that it is also the dual cone of $\overline{NE}(Z_{\cf})$, the
closure of $NE(Z_{\cf})$ in $A(Z_{\cf})$.


Given a divisor $D$ on $Z_{\cf}$, we shall denote by $[D]$ its class
in the Picard group of $Z_{\cf}$ and, also, its image into
$A(Z_{\cf})$. For a curve $C$ on $\gp^2$, $\tilde{C}$ (respectively,
$C^*$) will denote its strict (respectively, total) transform on the
surface $Z_{\cf}$ via the sequence of blow-ups given by $\cb_{\cf}$.
It is well known that the set $\{[L^*]\}\cup \{ [E_{q}^*]\}_{q\in
{\cal B}_{\cf}}$ is a $\gz$-basis (respectively, $\gr$-basis) of
$\pic(Z_{\cf})$ (respectively,  $A(Z_{\cf})$), where $L$ denotes a
general line of $\gp^2$.

Now assume that the foliation ${\cal F}$ has a rational first
integral. Then, we define the following divisor on the surface
$Z_{\cf}$:
$$
D_{\cf}:=d L^*-\sum_{q\in {\cal B}_{\cf}} r_q E_q^*,
$$
where $d$ is the degree of the curves in $\cp_{\cf}$ and $r_q$ the
multiplicity at $q$ of the strict transform  of a general curve of
$\cp_{\cf}$ on the surface that contains $q$.
Notice that the image on $A(Z_{\cf})$ of the
strict transform of a general curve of ${\cal P}_{\cf}$ coincides
with $[D_{\cf}]$ and, moreover, $D_{\cf}^2=0$ by Bézout Theorem.

Since $|D_{\cal F}|$ is a base-point-free complete linear system,
$NE(Z_{\cf})\cap [D_{\cf}]^\bot$ is the face of the cone
$NE(Z_{\cf})$ spanned by the images, in $A(Z_{\cal F})$, of those
curves on $Z_{\cf}$ contracted by the morphism $Z_{\cal
F}\rightarrow \gp |D_{\cal F}|$ determined by a basis of $|D_{\cal
F}|$. This will be a useful fact in this paper. In order to study
that face, we shall apply Cayley-Bacharach Theorem  \cite[
CB7]{eis}, which deals with residual schemes with respect to
complete intersections.


\begin{lem}\label{lem1}
If ${\cal F}$ is a foliation on $\gp^2$ with a rational first
integral, then the following equality, involving the projective
space of one dimensional quotients of
global sections of a sheaf of $\co_{\gp^2}$-modules, holds:
$$\gp H^0(\gp^2,{\pi_{\cf}}_*\co_{Z_{\cf}}(D_{\cf}))=\cp_{\cf}.$$
\end{lem}
\noindent {\it Proof}. Consider the ideal sheaf ${\cal J}$ defined
locally by the equations of the curves corresponding to two
general elements of the pencil ${\cal P}_{\cf}$ and its associated
zero dimensional scheme $\Gamma$.
${\pi_{\cf}}_*\co_{Z_{\cf}}(-\sum_{q\in {\cal B}_{\cf}} r_q
E_q^*)$ coincides with the integral closure of ${\cal J}$,
$\overline{\cal J}$, that is, the ideal sheaf such that in any
point  $p\in \gp^2$, the stalk ${\overline{\cal J}}_p$ is the
integral closure of the ideal ${\cal J}_p$ of $\co_{{\gp^2},p}$.
Applying Cayley-Bacharach Theorem to $\Gamma$ and the subschemes
defined by $\overline{\cal J}$ and ${\cal J}'=$
Ann$({\overline{\cal J}}/{\cal J})$, one gets:

$$
h^1(\gp^2, \overline{\cal J}(d))=h^0(\gp^2, {\cal
J}'(d-3))-h^0(\gp^2,  {\cal J}(d-3)),
$$
where $h^i$ means $\dim H^i$ $(0 \leq i \leq 1)$. The last term of
the above equality vanishes by Bézout Theorem and $\overline{\cal
J}(d)$ coincides with the sheaf
${\pi_{\cf}}_*\co_{Z_{\cf}}(D_{\cf})$. Moreover, ${\cal J}'$ is
nothing but the conductor sheaf of ${\cal J}$ (i.e., the sheaf that
satisfies that for all $p\in \gp^2$, the stalk of ${\cal J}'$ at $p$
is the common conductor ideal of the generic elements of ${\cal
J}_p$). Hence, $h^0(\gp^2, {\cal J}'(d-3))$ coincides with $p_g$,
the geometric genus of a general curve of ${\cal P}_{\cf}$. So, we
have the following equality:
$$h^1(Z_{\cf}, \co_{Z_{\cf}}(D_{\cf}))=h^1(\gp^2, {\pi_{\cf}}_*\co_{Z_{\cf}}(D_{\cf}))=p_g.$$
By Bertini's Theorem, the strict transform on $Z_{\cf}$ of any
general curve of ${\cal P}_{\cf}$ is smooth and so, its geometric
and arithmetic genus coincide. Therefore, using the Adjunction
Formula, we obtain:
$$h^1(Z_{\cf}, \co_{Z_{\cf}}(D_{\cf}))=1+(K_{Z_{\cf}}\cdot
D_{\cf})/2,$$ where $K_{Z_{\cf}}$ denotes a canonical divisor of
$Z_{\cf}$. Finally, the result follows by applying Riemann-Roch
Theorem to the divisor $D_{\cf}$.\findemo


\begin{teo}\label{teorema}

Let ${\cal F}$ be a foliation on $\gp^2$ with a rational first
integral. Then,

\begin{itemize}

\item[(a)] The image in $A(Z_{\cf})$ of a curve $C\hookrightarrow
Z_{\cf}$ belongs to $NE(Z_{\cf})\cap [D_{\cf}]^\bot$ if, and only
if, $C=D+E$ where $E$ is a sum (may be empty) of strict transforms
of non-dicritical exceptional divisors and, either $D=0$, or $D$ is
the strict transform on $Z_\cf$ of an invariant by $\cf$ curve.

\item[(b)] If $C$ is a curve on $Z_{\cf}$ that belongs to
$NE(Z_{\cf})\cap [D_{\cf}]^\bot$, then $C^2\leq 0$. Moreover,
$C^2=0$ if, and only if, $C$ is linearly equivalent to $rD_{\cf}$
for some positive rational number $r$.

\end{itemize}

\end{teo}

\noindent {\it Proof.} The morphism $h_{\cf}:Z_{\cf}\rightarrow
\gp^1$ defined by the sequence (\ref{seq2}) is induced by a linear
system $V \subseteq |D_{\cf}|$  such that the direct image by
${\pi_{\cf}}$ of rational functions induces a one-to-one
correspondence between $V$ and $\cp_{\cf}$ \cite[Th. II.7]{beau}.
Hence, by Lemma \ref{lem1}, $V=|D_{\cf}|$. Now, Clause $(a)$
follows from the fact that the integral curves contracted by
$h_{\cf}$ are the integral components of the strict transforms of
the curves belonging to the pencil $\cp_{\cf}$ and the strict
transforms of the non-dicritical exceptional divisors (see
\cite[Prop. 2.5.2.1]{julio} and \cite[Exer. 7.2]{casas}).

To show Clause $(b)$, consider an ample divisor $H$ on $Z_{\cf}$ and
the set $\Theta:=\{z\in A(Z_{\cf})\mid z^2>0 \;\; \mbox{and} \;\;
[H]\cdot z>0\}$.  $\Theta$ is contained in $NE(Z_{\cf})$ by
\cite[Cor. 1.21]{kollar} and, by Clause $(a)$,
 $[C]$ belongs to the boundary of $NE(Z_{\cf})$.
Hence, the first statement of Clause $(b)$ holds.

Finally, in order to prove the second statement of $(b)$, assume
that $C^2=0$. By \cite[ Rem. V.1.9.1]{hart}, there exists a basis of
$A(Z_{\cf})$ in terms of which the topological closure of $\Theta$,
$\overline{\Theta}$, gets the shape of a half-cone over an Euclidean
ball of dimension the cardinality of $\cb_{\cf}$. The strict
convexity of this Euclidean ball implies that the intersection of
the hyperplane $[D_{\cf}]^\bot$ with $\overline{\Theta}$ is just the
ray $\gr_{\geq 0}[D_{\cf}]$, where $\gr_{\geq 0}$ denotes the set of
non-negative real numbers. Therefore, $[C] \in \gr_{\geq
0}[D_{\cf}]$ by (a) and so, $[C]=r[D_{\cf}]$ for some
positive rational number $r$. \findemo\\

\begin{rem}
{\rm As a consequence of Theorem \ref{teorema}, next we give  two
conditions which, in case that one of them be satisfied, allow to
discard the existence of a rational first integral for a foliation
$\cf$ on $\gp^2$:
\begin{itemize}
\item[1. ] There exists an invariant curve,
$C$,  such that $\tilde{C}^2>0$. \item[2. ] There exist two
invariant curves, $C_1$ and $C_2$, such that $\tilde{C}_1^2=0$ and
$\tilde{C_1}\cdot \tilde{C_2}\not=0$.
\end{itemize}
}
\end{rem}


An essential concept for this paper is introduced in the following
\begin{de}
\label{condition} {\rm An {\it independent system of algebraic
solutions of  a foliation ${\cal F}$ on $\gp^2$}, which needs not to
have a rational first integral, is  a set  of  algebraic solutions
of ${\cal F}$, $S=\{C_1,C_2,\ldots,C_s\}$, where $s$ is the number
of dicritical exceptional divisors appearing in the minimal
resolution of ${\cal F}$, such that $\tilde{C}_i^2\leq 0$ ($1 \leq i
\leq s$) and the set of classes ${\cal A}_S:=\{[\tilde{C}_1],
[\tilde{C}_2],\ldots,[\tilde{C}_s]\}\cup \{[\tilde{E}_{q}]\}_{q\in
{\cal N}_{\cf}}\subseteq A(Z_{\cf})$ is $\gr$-linearly independent.}
\end{de}

\begin{rem}
{\rm Note that, when ${\cal F}$ has a rational first integral, the
existence of an independent system of algebraic solutions is an
equivalent fact to say that the face of the cone of curves of
$Z_{\cal F}$ given by $NE(Z_{\cal F})\cap [D_{\cal F}]^\bot$ has
codimension 1. }
\end{rem}

Assume now that a foliation ${\cal F}$ on $\gp^2$ admits an
independent system of algebraic solutions
$S=\{C_1,C_2,\ldots,C_s\}$. Set ${\cal
B}_{\cf}=\{q_1,q_2,\ldots,q_m\}$, ${\cal
N}_{\cf}=\{q_{i_1},q_{i_2},\ldots, q_{i_l}\}$ and stand $c_i:=(d_i,
-a_{i1}, \ldots, -a_{im})$, (respectively, $e_{q_{i_k}}: =(0,b_{k1},
\ldots, b_{km})$) for the coordinates of the classes of the strict
transforms on $Z_{\cf}$, $[\tilde{C}_i]$ (respectively,
$[\tilde{E}_{q_{i_k}}]$) of the curves $C_i$, $1 \leq i \leq s$,
(respectively, non-dicritical exceptional divisors $E_{q_{i_k}}$,
$1\leq k \leq l$) in
the basis of $A(Z_{\cf})$ given by $\{[L^*],[{E}^*_{q_1}],
[{E}^*_{q_2}],$ $\ldots,[{E}^*_{q_m}]\}$.

Notice that $d_i$ is the degree of $C_i$, $a_{ij}$  the multiplicity
of the strict transform of $C_i$ at $q_j$, and $b_{kj}$ equals 1 if
$j=i_k$, $-1$ if $q_j$ is proximate to $q_{i_k}$ and 0, otherwise.

Consider the divisor on $Z_{\cf}$:
\begin{equation}
\label{t} T_{{\cal F},S}:=\delta_0 L^*-\sum_{j=1}^m \delta_j
E_{q_j}^*,
\end{equation}
where $\delta_j:= \delta'_j /{\rm gcd}(\delta'_0, \delta'_1, \ldots,
\delta'_m)$, $\delta'_j$ being the absolute value of the determinant
of the matrix obtained by removing the $(j+1)$th column of the $m
\times (m+1)$-matrix defined by the rows $c_1, \ldots, c_s,
e_{q_{i_1}}, \ldots, e_{q_{i_l}}$.
Also, the set
$$
\Sigma({\cal F},S):=\{\lambda\in \gz_+
\mid h^0(\gp^2,{\pi_{\cf}}_*\co_{Z_{\cf}}(\lambda T_{{\cal F},S}))\geq
2\},
$$
where $h^0$ means $\dim H^0$ and $\gz_+$ is the set of positive
integers.

When the foliation ${\cal F}$ has a rational first integral, the
set ${\cal A}_S$
spans the hyperplane $[D_{\cf}]^\bot$ and, therefore,
$\sum_{i=0}^m \delta_i x_i=0$ is an equation for it, whenever
$(x_0,x_1,\ldots,x_m)$ are coordinates in $A(Z_{\cf})$ with
respect to the basis $\{[L^*]\}\cup \{ [{E}^*_q]\}_{q\in {\cal
B}_{\cf}}$. Hence, in this case, the divisor $D_{\cal F}$ is a
positive multiple of $T_{{\cal F},S}$. In fact, $[T_{{\cal F},S}]$
is the primitive element of the ray in $A(Z_{\cal F})$ spanned by
$[D_{\cal F}]$ in the sense that every divisor class belonging to
this ray is the product of $[T_{{\cal F},S}]$ by a positive
integer. Therefore the divisor $T_{{\cal F},S}$ does not depend on
the choice of the independent system of algebraic solutions $S$.

\begin{lem}\label{simplicial}
Consider a foliation ${\cal F}$ on $\gp^2$, with a rational first
integral, such that it admits an independent system of algebraic
solutions $S=\{C_j\}_{j=1}^s$.  Then,  $NE(Z_{\cal F})\cap
[D_{\cal F}]^\bot$ is a simplicial cone if the decomposition of
the class $[T_{{\cal F},S}]$ as a linear combination of the
elements in the set ${\cal A}_S$
contains every class in ${\cal A}_S$ and all its coefficients are
strictly positive.

\end{lem}
\noindent {\it Proof}. It follows from the fact that the convex
cone $NE(Z_{\cal F})\cap [D_{\cal F}]^\bot$ is spanned by the
classes in ${\cal A}_S$.
Indeed, by Theorem \ref{teorema}, it is enough to prove that, for
each algebraic solution $D\hookrightarrow \gp^2$ which does not
belong to $S$, the class $[\tilde{D}]$ can be written as a positive
linear combination of the above mentioned classes. But, since
$D_{\cf}\cdot \tilde{D}=0$ and $D_{\cf}$ is nef and a positive
multiple of $T_{\cf,S}$, one has that $\tilde{C_i}\cdot \tilde{D}=0$
for all $i=1,2,\ldots,s$ and $\tilde{E}_q\cdot \tilde{D}=0$ for each
$q\in {\cal N}_{\cf}$. Then, $[\tilde{D}]$ belongs to the subspace
of $A(Z_{\cf})$ orthogonal to ${\cal A}_S$ and, therefore, it must
be a positive multiple of $[T_{S,\cf}]$,
fact that proves the statement.\findemo\\

Now, we introduce some functions which will be instrumental in
stating our next result. Given a positive integer $k$, $\cd(k)$
will stand for the set of positive integers that divide $k$,
$\phi_{\cd(k)}$ (respectively, $-\phi_{\cd(k)}$) for the function
from the power set of $\cd(k)$, $\cp(\cd(k))$, to the rational
numbers, given by $\phi_{\cd(k)} (\Sigma) = 1-\sum_{\sigma \in
\Sigma} \frac{\sigma-1}{\sigma}$, (respectively, $-\phi_{\cd(k)}
(\Sigma) = (\sum_{\sigma \in \Sigma} \frac{\sigma-1}{\sigma})- 1$)
and   $w_{\cd(k)}: \gz \setminus\{0\} \rightarrow \gq_+$, $\gq_+$
being the set of positive rational numbers, for the assignment
defined by
\[
w_{\cd(k)}(a)=\left\{
\begin{array}
[c]{ll}
\min [\phi_{\cd(k)} (\cp(\cd(k)))\cap \gq_+] & \mbox{ if } a >0\\
\min [-\phi_{\cd(k)} (\cp(\cd(k)))\cap \gq_+] & \mbox{ if } a
<0\mbox{.}
\end{array}
\right.
\]

Notice that $w_{\cd(k)}$ is not defined neither for $a=0$ nor for
negative values when the set  $-\phi_{\cd(k)} (\cp(\cd(k)))\cap
\gq_+$ is empty.
\begin{teo}\label{gordo2}

Let ${\cal F}$ be a foliation on $\gp^2$ with a rational first
integral. Assume that ${\cal F}$ admits an independent system of
algebraic solutions $S=\{C_i\}_{i=1}^s$ and set
\begin{equation}
\label{dec} [T_{{\cal F},S}]=\sum_{i=1}^s \alpha_i [\tilde{C}_i]+
\sum_{q\in {\cal N}_{\cf}} \beta_q [\tilde{E}_{q}]
\end{equation} the
decomposition of $[T_{{\cal F},S}]$ as a linear combination of the
classes  in ${\cal A}_S$.
Then, the following properties hold:\\ 

(a) $D_{\cal F}=\alpha T_{{\cal F},S}$, where $\alpha$ is the
minimum of $\Sigma({\cal F},S)$.\\


(b) Assume that  the coefficients $\alpha_i$ $(1\leq i \leq s)$ and
$\beta_q$ ($q \in {\cal N}_{\cf}$) of the decomposition (\ref{dec})
are positive. Let $r$ be the minimum positive integer such that
$r\alpha_i, r\beta_q\in \gz$ for $i=1,2,\ldots,s$ and for $q\in
{\cal N}_{\cf}$, and let $k_0$ be the greatest common divisor of the
integers of the set $\{r T_{\cf,S}\cdot L^*\}\cup \{r T_{\cf,S}\cdot
E_q^*\}_{q\in {\cal B}_{\cf}}$. Then

\begin{equation}\label{equat}
\alpha \leq \Delta_{\cf}:=\frac{\deg({\cal F})+2-\sum_{i=1}^s
\deg(C_i)}{w_{\cd(k_0)}(\deg({\cal F})+2-\sum_{i=1}^s \deg(C_i))
\sum_{i=1}^s \alpha_i \deg(C_i)},
\end{equation}
where $\deg({\cal F})$ denotes the degree of the foliation ${\cal
F}$ and  $\deg(C_i)$  the one of the curve $C_i$ ($1 \leq i \leq
s$).



\end{teo}

\noindent {\it Proof}.  Let $\mu$ be the positive integer such
that $D_{\cf}=\mu T_{{\cal F},S}$.

In order to prove $(a)$ we shall reason by contradiction assuming
that $\alpha<\mu$. Taking into account that $D_{\cf}\cdot T_{{\cal
F},S}=0$ and $D_{\cf}$ is nef, and applying Theorem \ref{teorema},
it is deduced that all the elements of the linear system $\gp
H^0(\gp^2,{\pi_{\cf}}_*\co_{Z_{\cf}}(\alpha T_{{\cal F},S}))$ are
invariant curves. So, its integral components must be also integral
components of the fibers of the pencil $\cp_{\cf}$. It is clear that
there exist infinitely many integral components of elements in $\gp
H^0(\gp^2,{\pi_{\cf}}_*\co_{Z_{\cf}}(\alpha T_{{\cal F},S}))$ whose
strict transforms have the same class in the Picard group of
$Z_{\cf}$. Finally, since $\cp_{\cf}$ is an irreducible pencil, the
unique class in $\pic(Z_{\ck})$ corresponding to an infinite set of
integral curves is that of the general fibers of $\cp_{\cf}$, which
is a contradiction, because the degree of these general fibers is
larger than the degree of the curves in $\gp
H^0(\gp^2,{\pi_{\cf}}_*\co_{Z_{\cf}}(\alpha T_{{\cal F},S}))$.

Next, we shall prove $(b)$. By Lemma \ref{simplicial}, $NE(Z_{\cal
F})\cap [D_{\cal F}]^\bot$ is the simplicial convex cone spanned by
the classes  in $\{[\tilde{C}_i]\}_{1 \leq i \leq s}  \cup
\{[\tilde{E}_{q}]\}_{q\in {\cal N}_{\cf}}$. Firstly, we shall show
that $\sum_{i=1}^s \mu \alpha_i C_i$ is the unique curve of ${\cal
P}_{\cal F}$ containing some curve of $S$ as a component.

To do it, assume that $Q$ is a fiber of ${\cal P}_{\cal F}$
satisfying the mentioned condition. Then, by Clause $(a)$ of Theorem
\ref{teorema}, there exists a sum of strict transforms of
non-dicritical exceptional divisors, $E$, such that $\tilde{Q}+E$ is
a divisor linearly equivalent to $D_{\cf}$. Taking the decomposition
of $Q$ as a sum of integral components, one gets:
\begin{equation}\label{qq}
\tilde{Q}+E=\sum_{i=1}^s a_i \tilde{C}_i+\sum_{j=1}^t b_j
\tilde{D}_j+\sum_{q\in {\cal N}_{\cf}} c_q \tilde{E}_q,
\end{equation}
where all coefficients are non-negative integers, some $a_i$ is
positive and the elements in $\{D_j\}_{1\leq j \leq t}$ are integral
curves on $\gp^2$ which are not in $S$. If some $\tilde{D}_j$ had
negative self-intersection, its image in $A(Z_{\cf})$ would span an
extremal ray of the cone $NE(Z_{\cal F})\cap [D_{\cal F}]^\bot$,
contradicting the fact that this cone is simplicial. Hence, by
Theorem \ref{teorema}, $\tilde{D}_j^2=0$ for all $j=1,2,\ldots ,t$
and all the classes $[\tilde{D}_j]$ belong to the ray spanned by
$[D_{\cf}]$. As a consequence, one has the following decomposition:
$$[\tilde{Q}]+[E]=\sum_{i=1}^s a_i [\tilde{C}_i]+\sum_{j=1}^t b'_j
\left(\sum_{i=1}^s \alpha_i [\tilde{C}_i]+ \sum_{q\in {\cal
N}_{\cf}} \beta_q [\tilde{E}_{q}]\right)+\sum_{q\in {\cal
N}_{\cf}} c_q [\tilde{E}_q]=$$
$$=\sum_{i=1}^s \left(a_i+ \alpha_i\sum_{j=1}^t b_j'\right)[\tilde{C}_i]+\sum_{q\in {\cal N}_{\cf}} \left(c_q+ \beta_q\sum_{j=1}^t b_j'\right)[\tilde{E}_q],$$
where, for each $j=1,2,\ldots,t$, $b_j'$ is  $b_j$ times a positive
integer. As $[\tilde{Q}]+[E]=[D_{\cf}]=\mu ( \sum_{i=1}^s \alpha_i
[\tilde{C}_i]+ \sum_{q\in {\cal N}_{\cf}} \beta_q [\tilde{E}_{q}]
)$, $a_i=(\mu-\sum_{j=1}^t b_j') \alpha_i$ for all $i=1,2,\ldots,s$
and $c_q=(\mu-\sum_{j=1}^t b_j') \beta_q$ for all $q\in {\cal
N}_{\cf}$. Then, since some $a_i$ does not vanish and all the
rational numbers $\alpha_i$ are different from zero, it follows that
$a_i\not=0$ for all $i=1,2,\ldots,s$. Now, from the inequality
$$
s-1\geq \sum_{Q} (n_Q-1),
$$
where the sum is taken over the set of fibers $Q$ of ${\cal
P}_{\cf}$ and $n_Q$ denotes the number of distinct integral
components of $Q$ \cite{kaliman}, we deduce that $Q$ cannot have
integral components different from those in $S$ and, therefore, we
must take, in the equality (\ref{qq}), $b_j=0$ for each
$j=1,2,\ldots,t$, $a_i=\mu \alpha_i$ ($1\leq i\leq s$) and $c_q=\mu
\beta_q$ for all $q\in {\cal N}_{\cf}$. As a consequence, we have
proved the mentioned property of the curve $\sum_{i=1}^s \mu\alpha_i
C_i$ and, moreover, that the integer $r$ defined in the statement
divides $\mu$.

Secondly, from the above paragraph it can be deduced that the
non-integral fibers of $\cp_{\cf}$ different from $\sum_{i=1}^s
\mu\alpha_i C_i$ must have the form $nD$, where $n>1$ is a
positive integer, $D$ is an integral curve and $\tilde{D}^2=0$.
Next, we bound such an integer $n$.

Consider the effective divisor of $Z_{\cf}$, $R:=\sum_{i=1}^s
r\alpha_i \tilde{C}_i+ \sum_{q\in {\cal N}_{\cf}} r\beta_q
\tilde{E}_{q}$. $[\tilde{D}]$ and $[R]$ span the same ray of
$A(Z_{\cf})$ by Clause $(b)$ of Theorem \ref{teorema} and so,
there exist two relatively prime positive integers $a,b$  such
that $[\tilde{D}]=\frac{a}{b}[R]$. On the one hand,
$\frac{\mu}{r}[R]=[D_{\cf}]=\frac{na}{b}[R]$ and, since
$\frac{\mu}{r}\in \gz$, one gets that $b$ divides $n$. On the
other hand, $b\tilde{D}$ and $aR$ are two linearly equivalent
effective divisors without common components. Then, the linear
system of $\gp^2$ given by
$H^0(\gp^2,{\pi_{\cf}}_*\co_{Z_{\cf}}(b\tilde{D}))$ has no fixed
components and all its members are invariant curves by $\cf$, by
Theorem \ref{teorema}. Therefore, either it has irreducible
general members (in which case it coincides with $\cp_{\cf}$), or
it is composite with an irreducible pencil (that must be
$\cp_{\cf}$). Hence, there exists a positive integer $\delta$ such
that $b[\tilde{D}]=\delta[D_{\cf}]$. But, since
$[D_{\cf}]=n[\tilde{D}]$, one has that $n$ divides (and then it is
equal to) $b$. As a consequence,  $n[\tilde{D}]=a[R]$ and, since
$\gcd(a,n)=1$, $n$ divides the coordinates of $[R]$ with respect
to the basis $\{[L^*]\}\cup\{E_q^*\}_{q\in {\cal B}_{\cf}}$, that
is, $n$ must divide the integer $k_0$ defined in the statement.

Finally, (b) follows from the  formula that asserts that if a plane
foliation ${\cal F}$ has a rational first integral, with associated
pencil of degree $d$, then $2d- \deg {\cal F} -2=\sum (e_Q-1)
\deg(Q)$, where the sum is taken over the set of integral components
$Q$ of the curves in $\cp_{\cf}$ and $e_Q$ denotes the multiplicity
of $Q$ as component of the correspondent fiber. Indeed, by applying
this formula one gets:
$$\mu = \frac{\deg({\cal F})+2-\sum_{i=1}^s
\deg(C_i)}{(1-\sum_{Q\in \Upsilon} \frac{e_Q-1}{e_Q})\sum_{i=1}^s
\alpha_i \deg(C_i)},$$ where $\Upsilon$ denotes the set of
integral curves $D$ which do not belong to $S$ and such that $nD$
is a member of the pencil $\cp_{\cf}$ for some positive integer $n
>1$. Notice that,  by $(a)$, $h^0(Z_{\cal F}, \co_{Z_{\cal
F}}(D))=1$ and so the function from $\Upsilon$ to the positive
integers given by $D\mapsto n$ is injective. Now, the bound
$\Delta_{\cf}$ of the statement is clear, since the integers
$e_Q$ (if they exist) are divisors of $k_0$.\findemo\\

\begin{rem}\label{well}
{\rm For an arbitrary foliation admitting an independent system of
algebraic solutions $S$ as above, the value $\Delta_{\cf}$ given in
Theorem \ref{gordo2} could be less than 1 or even not computable
because $w_{\cd(k_0)}(\deg({\cal F})+2-\sum_{i=1}^s \deg(C_i))$
could not be  defined. In those cases, $\cf$ has no rational first
integral and we say that $\Delta_{\cf}$  \it{is not well defined}. }
\end{rem}







Stands $K_{Z_{\cf}}$ for a canonical divisor on $Z_{\cf}$. The
next result follows from Bertini's Theorem and the Adjunction
Formula, and it shows that the condition
 $K_{Z_{\cf}}\cdot T_{\cf,S}<0$ makes easy to check
whether  $\cf$ has or not a rational first integral, and to
compute it (using Lemma \ref{lem1}).

\begin{pro}\label{memo}

Let $\cf$ be a foliation on $\gp^2$ admitting an independent
system of algebraic solutions $S$. Assume that $K_{Z_{\cf}}\cdot
T_{\cf,S}<0$ and $\cf$ has a rational first integral. Then, the
general elements of the pencil ${\cal P}_{\cal F}$ are rational
curves and $D_{\cal F}=T_{\cf,S}$.

\end{pro}

\section{Computation of rational first integrals. Algorithms and
examples}\label{TRES}

With the above notations, recall that, for each foliation $\cf$ on
$\gp^2$ with a rational first integral, the divisor $D_{\cf}$
satisfies $D_{\cf}^2=0$ and $D_{\cf}\cdot \tilde{E}_{q}=0$
(respectively, $D_{\cf}\cdot \tilde{E}_{q}>0$) for all $q \in
{\cal N}_{\cf}$ (respectively, $q \in {\cal B}_{\cf}\setminus
{\cal N}_{\cf})$ \cite[Exer. 7.2]{casas}. These facts and Lemma
\ref{lem1} support the following decision algorithm for the
problem of deciding whether an arbitrary foliation $\cf$ has a
rational first integral of a fixed degree $d$. In fact, it allows
to compute it in the affirmative case.

This is an alternative algorithm to the particular case of that
given by Prelle and Singer \cite{pr-si} (and implemented by Man
\cite{man1}) to compute first integrals of foliations. Note that our
algorithm only involves integer arithmetic and resolution of systems
of linear equations (Step 3), and that we do not need to use
Groebner bases.

 Let $\cf$ be a foliation on $\gp^2$ and $d$ a positive
integer.\\
\begin{alg}
\label{alg1} $\;$\newline
 \noindent {\rm {\it Input:} $d$, a projective 1-form $\Omega$
defining ${\cal F}$, ${\cal B}_{\cf}$ and ${\cal N}_{\cf}$.

\noindent {\it Output:} Either a rational first integral of degree
$d$, or ``0'' if there is no such first integral.

\begin{itemize}

\item[1.] Compute the finite set $\Gamma$ of divisors
$D=dL^*-\sum_{q\in {\cal B}_{\cf}} e_q E_q^*$ such that

\begin{itemize}
\item[(a)] $D^2=0$,

\item[(b)] $D\cdot \tilde{E}_q=0$ for all $q\in {\cal N}_{\cf}$
and

\item[(c)] $D\cdot \tilde{E}_q> 0$ for all $q\in {\cal
B}_{\cf}\setminus {\cal N}_{\cf}$.

\end{itemize}

\item[2.] Pick $D\in \Gamma$.

\item[3.] If  the dimension of the $k$-vector space
$H^0(\gp^2,{\pi_{\cf}}_*\co_{Z_{\cf}}(D))$ is 2, then take a basis
$\{F,G\}$ and check the condition $d(F/G)\wedge \Omega=0$. If it
is satisfied, then return $F/G$.

\item[4.] Set $\Gamma:=\Gamma\setminus \{D\}$.

\item[5.] Repeat the steps 2,  3 and 4 while the set $\Gamma$ is not
empty.

\item[6.] Return ``0''.

\end{itemize}

}
\end{alg}

Now, and again supported in results of the previous section, we
give an algorithm to decide, under certain conditions, whether an
{\it arbitrary foliation} ${\cal F}$ on $\gp^2$ has a rational
first integral (of arbitrary degree). As above, the algorithm
computes it in the affirmative case.

Firstly we state the needs in order that  the algorithm works. Let
${\cal F}$ be a foliation on $\gp^2$ that admits an independent
system of algebraic solutions $S$ which satisfies {\it at least
one} of the following conditions:

\begin{itemize}

\item[(1)] $T_{\cf,S}^2\not=0$.

\item[(2)] The decomposition of the class $[T_{\cf,S}]$ as a
linear combination of those in the set ${\cal A}_S$, given in
Definition \ref{condition}, contains all the classes in ${\cal A}_S$
with positive coefficients. In this case, it will be useful the
value $\Delta_{\cf}$ given in Theorem \ref{gordo2}.


\item[(3)] The set $\Sigma({\cal F},S)$ is not empty.

\end{itemize}

\begin{alg}\label{alg2}
$\;$ \newline {\rm \noindent {\it Input:} A projective 1-form
$\Omega$ defining ${\cal F}$, ${\cal B}_{\cf}$, ${\cal N}_{\cf}$
and an independent system of algebraic solutions $S$ satisfying at
least one of the above conditions $(1)$, $(2)$ and $(3)$.

\noindent {\it Output:} A rational first integral for ${\cal F}$,
or ``0'' if there is no such first integral.

\begin{itemize}

\item[1.] If  $(1)$ holds, then return
``0''.

\item[2.] If Condition $(2)$ is satisfied and either
$\Delta_{\cf}$ is not  well defined (see Remark \ref{well}) or
$h^0(\gp^2,{\pi_{\cf}}_*\co_{Z_{\cf}}(\lambda T_{{\cal F},S}))\leq
1$ for all positive integer $\lambda \leq \Delta_{\cf}$, then
return ``0''.

\item[3.] Let $\alpha$ be the minimum of the set $\Sigma(\cf,S)$.

\item[4.] If $h^0(\gp^2,{\pi_{\cf}}_*\co_{Z_{\cf}}(\alpha T_{{\cal
F},S}))>2$,  then return ``0''.

\item[5.] Take a basis $\{F,G\}$ of
$H^0(\gp^2,{\pi_{\cf}}_*\co_{Z_{\cf}}(\alpha T_{{\cal F},S}))$ and
check the equality $d(F/G)\wedge \Omega=0$. If it is satisfied,
then return $F/G$. Else, return ``0''.

\end{itemize}

}

\end{alg}

\begin{exa}\label{e1}
{\rm
 Let ${\cal F}$ be the foliation on the projective plane over the complex numbers defined by the
 projective 1-form $\Omega=AdX+BdY+CdZ$, where
$$A=X^3Y+4Y^4+2X^3Z-X^2Y Z-4X^2 Z^2-X Y Z^2+2X Z^3+YZ^3,$$
$$B=-X^4-4X Y^3+3X^3Z+4Y^3Z-3X^2Z^2+X Z^3,$$
$$C= -2X^4-2X^3Y-4Y^4+4X^3Z+4X^2YZ-2X^2Z^2-2X Y Z^2.$$
Resolving $\cf$, we have computed the configuration $\ck_{\cf}$,
which coincides with the configuration ${\cal
B}_{\cf}=\{q_j\}_{j=1}^{10}$ of dicritical points. The proximity
graph is given in Figure 1. ${\cal N}_{\cf} =
\{q_1,q_2,q_3,q_4,q_5,q_7,q_8,q_9\}$ and $S=\{C_1,C_2\}$ is an
independent system of algebraic solutions, where $C_1$
(respectively, $C_2$) is the line (respectively, conic) given by the
equation $X-Z=0$ (respectively, $(8i-1)X^2+4i X Y+8Y^2+(2-8i)X Z-4iY
Z-Z^2=0$). The divisor $T_{\cf,S}$ is
$4L^*-2E_1^*-2E_2^*-\sum_{j=3}^{10}E_j^*$. Clearly, Condition (1)
above is not satisfied. (2) does not hold either, because
$$
[T_{\cf,S}]=4[\tilde{C}_1]+2[\tilde{E}_1]+4[\tilde{E}_2]+3[\tilde{E}_3]+2[\tilde{E}_4]
+[\tilde{E}_5]+3[\tilde{E}_7]+2[\tilde{E}_8]+[\tilde{E}_9].
$$
The space of global sections $H^0(\gp^2,{\pi_{\cf}}_*\co_{Z_{\cf}}(
T_{{\cal F},S}))$ has dimension $2$ and is spanned by
$F=X^2Z^2-2X^3Z+X^4+XYZ^2-2X^2YZ+X^3Y+Y^4$ and $G=(X-Z)^4$.
Therefore, Condition (3) happens and if ${\cal F}$ admits rational
first integrals, one of them must be $R:= F/G$.
The equality $\Omega \wedge dR=0$ shows that $R$ is, in fact, a
rational first integral of ${\cal F}$.
 }
\end{exa}

\begin{figure}[ht]\label{fig}
\setlength{\unitlength}{1mm}
\begin{center}
\begin{picture}(100,60)

\put(50,0){\circle*{2}} \put(50,0){\line(0,1){10}}
\put(55,0){$q_1$}

\put(50,10){\circle*{2}}\put(50,10){\line(1,1){10}}
\put(50,10){\line(-1,1){10}} \put(55,10){$q_2$}

\put(40,20){\circle*{2}}\put(40,20){\line(-1,1){10}}
\put(45,20){$q_3$}

\put(30,30){\circle*{2}}\put(30,30){\line(-1,1){10}}
\put(35,30){$q_4$}

\put(20,40){\circle*{2}}\put(20,40){\line(-1,1){10}}
\put(25,40){$q_5$}

\put(10,50){\circle*{2}} \put(15,50){$q_6$}

\put(60,20){\circle*{2}}\put(60,20){\line(1,1){10}}
\put(65,20){$q_7$}

\put(70,30){\circle*{2}}\put(70,30){\line(1,1){10}}
\put(75,30){$q_8$}

\put(80,40){\circle*{2}}\put(80,40){\line(1,1){10}}
\put(85,40){$q_9$}\put(90,50){\circle*{2}} \put(95,50){$q_{10}$}

\end{picture}
\end{center}
\caption{The proximity graph  of ${\cal B}_{\cf}$ in Example
\ref{e1}}
\end{figure}

\begin{rem}
{\rm Assume that $\cf$ has a rational first integral and admits an
independent system of algebraic solutions. Set $k(T_{\cf,S}) =
{\rm tr. deg.} (\oplus_{n \geq 0} H^0( {\cal O}_{Z_{\cf}} (n
T_{\cf,S}))) -1$ the $T_{\cf,S}$-dimension of $Z_{\cf}$.
$k(T_{\cf,S}) \leq 2$ and, by using results in \cite{zar} and
\cite{cut}, it can be proved that  Condition (3) given before
Algorithm \ref{alg2} only happens when $k(T_{\cf,S}) = 1$.}

\end{rem}

To decide whether an independent system of algebraic solutions
satisfies one of the above mentioned conditions (1) or (2) is very
simple, but to check Condition (3) may be more difficult. However,
when $K_{Z_{\cf}}\cdot T_{\cf,S}<0$, we should not be concerned
about these conditions since, by Proposition \ref{memo}, there is
no need to take all the steps in Algorithm \ref{alg2}. Indeed, it
suffices to check whether
$h^0(\gp^2,{\pi_{\cf}}_*\co_{Z_{\cf}}(T_{{\cal F},S}))=2$; in the
affirmative case, we shall go to Step 5 and, otherwise, $\cf$ has
no rational first integral. Next, we shall show an enlightening
example.

\begin{exa}\label{ex4}
{\rm

Set ${\cal F}$  the foliation on the complex projective plane
defined by the 1-form
$$\Omega = 2YZ^5dX + (-7Y^{5} Z - 3X Z^5+Y Z^5)dY + (7Y^{6} + X Y Z^4-Y^2
Z^4)dZ.$$ The configuration $\ck_{\cf}$ coincides with the one  of
dicritical points ${\cal B}_{\cf}$. It has 13 points and its
proximity graph is given in Figure 2.  The dicritical exceptional
divisors are $E_{q_3}$ and $E_{q_{13}}$. The  set $S$  given by the
lines with equations $Y=0$ and $Z=0$ is an independent system of
algebraic solutions. Its associated divisor $T_{\cf,S}$ is
$$
10L^*-2E_{q_1}^*-E_{q_2}^*-E_{q_3}^*-8E_{q_4}^*-2\sum_{i=5}^{11}
E_{q_i}^*-E_{q_{12}}^*-E_{q_{13}}^*.
$$
Now, by Proposition \ref{memo}, if ${\cal F}$ has a rational first
integral, then the pencil $\cp_{\cf}$ is $\gp
H^0(\gp^2,{\pi_{\cf}}_*\co_{Z_{\cf}}(T_{\cf,S}))$. A basis of this
projective space is given by $F_1=Y^{10}-2XY^5Z^4+2Y^6Z^4+X^2Z^8-2X
Y Z^8+Y^2 Z^8$ and $F_2=Y^3 Z^7$. Finally, the equality
$d(F_1/F_2)\wedge \Omega =0$ shows that $F_1/F_2$ is a rational
first integral of $\cf$.

}
\end{exa}

\begin{figure}[hbt]\label{fig}
\setlength{\unitlength}{1mm}
\begin{center}
\begin{picture}(40,45)

\put(0,0){\circle*{2}} \put(0,0){\line(0,1){5}} \put(5,0){$q_1$}

\put(0,5){\circle*{2}}\put(0,5){\line(0,1){5}} \put(5,5){$q_2$}

\put(0,10){\circle*{2}} \put(5,10){$q_3$}

\put(30,0){\circle*{2}} \put(30,0){\line(0,1){5}} \put(35,0){$q_4$}

\put(30,5){\circle*{2}} \put(30,5){\line(0,1){5}} \put(35,5){$q_5$}

\put(30,10){\circle*{2}} \put(30,10){\line(0,1){5}}
\put(35,10){$q_6$}

\put(30,15){\circle*{2}} \put(30,15){\line(0,1){5}}
\put(35,15){$q_7$}

\put(30,20){\circle*{2}} \put(30,20){\line(0,1){5}}
\put(35,20){$q_8$}

\put(30,25){\circle*{2}} \put(30,25){\line(0,1){5}}
\put(35,25){$q_9$}

\put(30,30){\circle*{2}} \put(30,30){\line(0,1){5}}
\put(35,30){$q_{10}$}

\put(30,35){\circle*{2}} \put(30,35){\line(0,1){5}}
\put(35,35){$q_{11}$}

\put(30,40){\circle*{2}} \put(30,40){\line(0,1){5}}
\put(35,40){$q_{12}$}

\put(30,45){\circle*{2}} \put(35,45){$q_{13}$}

\qbezier[30](30,0)(15,10)(30,20)

\qbezier[20](30,35)(20,40)(30,45)

\qbezier[20](0,0)(-10,5)(0,10)

\end{picture}
\end{center}
\caption{The proximity graph of ${\cal B}_{\cf}$ in Example
\ref{ex4}}
\end{figure}

The following example shows the existence of foliations on $\gp^2$
with a rational first integral which do not admit an independent
system of algebraic solutions.

\begin{exa}
{\rm Consider the foliation $\cf$ on the projective plane over the
complex numbers given by the projective 1-form that defines the
derivation on
the following rational function:
$$
\dfrac{XZ^2+3YZ^2-Y^3}{YZ^2+3XZ^2-X^3}.
$$
The configuration of dicritical points ${\cal B}_{\cf}$ has 9
points, all in $\gp^2$, and therefore, ${\cal N}_{\cf}=\emptyset$.
Moreover, the divisor $D_{\cf}$ is given by $3L^*-\sum_{q\in {\cal
B}_{\cf}} E_q^*$. By Theorem \ref{teorema}, the intersection product
of the strict transform on $Z_{\cf}$ of whichever invariant curve of
$\cf$ times $D_{\cal F}$ vanishes, and so
the non-general algebraic solutions are among the lines passing
through three points in ${\cal B}_{\cf}$ and the irreducible
conics passing through six points in ${\cal B}_{\cf}$. Simple
computations show that these are, exactly, the 8 curves given by
the equations:

$$
X-Y=0; \; X+Y=0; \; 2X+(\sqrt{5}+3)Y=0; \; -2X+(\sqrt{5}-3)Y=0;
$$
$$
X^2-X Y+Y^2-4Z^2=0; \; X^2+X Y+Y^2-2Z^2=0;
$$
$$
2X^2+(\sqrt{5}-3)X Y-(3 \sqrt{5}-7)Y^2+(8 \sqrt{5}-24)Z^2=0;
$$
$$
-2X^2+(\sqrt{5}+3)X Y-(3 \sqrt{5}+7)Y^2+(8 \sqrt{5}+24)Z^2=0.
$$
Hence, $\cf$ does not admit an independent system of algebraic
solutions.

}
\end{exa}

The following result
will help to state an algorithm, for foliations $\cf$ whose cone
$NE(Z_{\cf})$ is polyhedral, that either computes an independent
system of algebraic solutions or discards that $\cf$ has a
rational first integral. In the sequel, for each subset $W$ of
$A(Z_{\cf})$, $\con(W)$ will denote the convex cone of
$A(Z_{\cf})$ spanned by $W$.

\begin{pro}\label{nueva}
Let $\cf$ be a foliation on $\gp^2$ having a rational first integral
and such that $NE(Z_{\cf})$ is polyhedral. Let $G$ be a non-empty
finite set of integral curves on $\gp^2$ such that $x^2\geq 0$ for
each element $x$ in the dual cone $\con(W)^\vee$, $W$ being the
following subset of $A(Z_{\cf})$: $W=\{[\tilde{Q}]\mid Q\in G\}\cup
\{[\tilde{E}_q]\}_{q\in {\cal B}_{\cf}}$. Then, $\cal F$ admits an
independent system of algebraic solutions $S$ such that $S\subseteq
G$.
\end{pro}
\noindent {\it Proof}. The conditions of the statement imply
$\con(W)^\vee \subseteq \overline{\Theta}$, where
$\overline{\Theta}$ is the topological closure  of the set
$\Theta=\{x\in A(Z_{\cf})\mid x^2> 0 \mbox{ and } [H]\cdot x>0\}$
given in the proof of  Theorem \ref{teorema}. Thus,
$\overline{\Theta}^\vee \subseteq (\con(W)^\vee)^\vee=\con(W)$,
where the equality is due to the fact that $\con(W)$ is closed.
Now, $[D_{\cf}]\in P(Z_{\cf})\subseteq \overline{\Theta}^\vee$ by
\cite[Cor. 1.21]{kollar}, and so $[D_{\cf}]\in \con(W)$. Note that
the cones $\con(W)$ and $NE(Z_{\cf})$ have the same dimension and,
as $[D_{\cf}]$  is in  the boundary of $NE(Z_{\cf})$, it also
belongs to the boundary of $\con(W)$.
Moreover $[D_{\cf}]\in \con(W)^\vee$. As a consequence, $\con(W)\cap
[D_{\cf}]^\bot$ is a face of $\con(W)$ which, in addition, contains
the class $[D_{\cf}]$.

Let $R$ be the maximal proper face of $\con(W)$ containing
$\con(W)\cap [D_{\cf}]^\bot$. Since $\con(W)$ is polyhedral, there
exists $y\in \con(W)^\vee\setminus \{0\}$ such that $R=\con(W)\cap
y^\bot$. Note that $y^2\geq 0$ and $y\in [D_{\cf}]^\bot$.



Recalling the shape of $\overline{\Theta}$ (see the proof of
Theorem \ref{teorema}), it is clear that the hyperplane
$[D_{\cf}]^\bot$ is tangent to the boundary of the half-cone
$\overline{\Theta}$. Thus, for each $x\in [D_{\cf}]^\bot \setminus
\{0\}$, $x^2\leq 0$ and the following equivalence holds:
\begin{equation}\label{equiv}
\mbox{$x^2<0$ if, and only if, $x$ does not belong to the line
$\gr[D_{\cf}]$.}
\end{equation}
Therefore $y$ belongs to the ray $\gr_{\geq 0}[D_{\cf}]$ and so the
equality $R=\con(W)\cap [D_{\cf}]^\bot$ is satisfied, which
concludes the proof by taking into account Clause (a) of Theorem
\ref{teorema} and that $R$ is a face of codimension one.\findemo

\begin{cor}
\label{poli} Let $\cf$ be a foliation on $\gp^2$ having a rational
first integral and such that $NE(Z_{\cf})$ is polyhedral. Then,
${\cal F}$ admits an independent system of algebraic solutions
$S$. Moreover, $S$ can be taken such that $\tilde{C}^2<0$ for all
$C\in S$.
\end{cor}
\noindent {\it Proof}. The cone $NE(Z_{\cf})$ is strongly convex
(by Kleiman's ampleness criterion \cite{kle}) and closed, so it is
spanned by its extremal rays. Thus, setting $G$ the set of
integral curves on $\gp^2$ whose strict transforms on $Z_{\cf}$
give rise to generators of extremal rays of $NE(Z_{\cf})$ and
applying Proposition \ref{nueva}, we prove our first statement.
The second one follows simply by taking into account that, since
the cardinality of ${\cal B}_{\cf}$ is assumed to be larger than
1, the extremal rays of $NE(Z_{\cf})$ are just those spanned by
the classes of the integral curves on $Z_{\cf}$ with negative
self-intersection (due to the polyhedrality of $NE(Z_{\cf})$).
\findemo\\


Now, we state the announced algorithm, where $\cf$ is a foliation
on $\gp^2$ such that the cone $NE(Z_{\cf})$ is polyhedral.

\begin{alg}\label{alg3}
$\;$ \newline \noindent {\it Input:} {\rm A projective 1-form
$\Omega$ defining $\cf$, ${\cal B}_{\cf}$ and ${\cal N}_{\cf}$.

\noindent {\it Output:} Either ``0'' (which implies that $\cal F$
has no rational first integral) or an independent system of
algebraic solutions.

\begin{itemize}

\item[1.] Define $V:=\con(\{[\tilde{E}_q]\}_{q\in {\cal
B}_{\cf}})$, $G:=\emptyset$ and let $\Gamma$ be the set of
divisors $C=dL^*-\sum_{q\in {\cal B}_{\cf}} e_q E_q^*$ satisfying
the following conditions:

\begin{itemize}

\item[(a)]  $d>0$ and $0 \leq e_q\leq d$ for all $q\in {\cal
B}_{\cf}$.

\item[(b)] $C\cdot \tilde{E}_q\geq 0$ for all $q\in {\cal
B}_{\cf}$.

\item[(c)] Either $C^2=K_{Z_{\cf}}\cdot C=-1$, or $C^2<0$ and
$K_{Z_{\cf}}\cdot C\geq 0$.

\end{itemize}

\item[2.] Pick $D\in \Gamma$ such that $D\cdot L^*$ is minimal.

\item[3.] If $D$ satisfies the conditions

\begin{itemize}

\item[(a)] $[D]\not\in V$,

\item[(b)] $h^0(\gp^2,{\pi_{\cf}}_*\co_{Z_{\cf}}(D))=1$,

\item[(c)] $[D]=[\tilde{Q}]$, where $Q$ is the divisor of zeros of a global
section of ${\pi_{\cf}}_*\co_{Z_{\cf}}(D)$,

\end{itemize}

then set $V:=\con(V\cup \{[D]\})$. If, in addition, $Q$ is an
invariant curve of $\cal{F}$, no curve in $G$ is a component of
$Q$ and $\{[\tilde{R}]\mid R\in G \}\cup \{[D]\}\cup
\{[\tilde{E}_q]\}_{q\in {\cal N}_{\cf}}$ is a $\gr$-linearly
independent system of $A(Z_{\cf})$, then set $G:=G\cup \{Q\}$.

\item[4.] Let $\Gamma:=\Gamma\setminus \{D\}$.

\item[5.] Repeat the steps 2, 3 and 4 while the following two
conditions are satisfied:

\begin{itemize}

\item[(a)] $\card(G)<\card({\cal B}_{\cf}\setminus {\cal
N}_{\cf})$, where $\card$ stands for cardinality.

\item[(b)] There exists $x\in V^\vee$ such that $x^2<0$.

\end{itemize}

\item[6.] If $\card(G)<\card({\cal B}_{\cf}\setminus {\cal
N}_{\cf})$, then return ``0". Else, return $G$.
\end{itemize}

}
\end{alg}

\noindent {\it Explanation.
} \label{explanation} This algorithm computes a strictly increasing
sequence of convex cones $V_0 \subset V_1\subset \cdots$ such that
$V_0=\con(\{[\tilde{E}_q]\}_{q\in {\cal B}_{\cf}})$ and
$V_i=\con(\{V_{i-1} \cup [\tilde{Q}_i]\})$ for $i\geq 1$, where $
Q_1, Q_2,\dots$ are curves on $\gp^2$ satisfying the following
conditions:
\begin{itemize}
\item[1)] Either $\tilde{Q}_i^2=K_{Z_{\cf}}\cdot \tilde{Q}_i=-1$,
or $\tilde{Q}_i^2<0$ and $K_{Z_{\cf}}\cdot \tilde{Q}_i\geq 0$.

\item[2)] $h^0(\gp^2,{\pi_{\cf}}_*\co_{Z_{\cf}}(\tilde{Q}_i))=1$.

\end{itemize}

Notice that, since the cone $NE(Z_{\cal F})$ is polyhedral, it has
a finite number of extremal rays and, moreover, they are generated
by the classes of the integral curves of $Z_{\cal F}$ having
negative self-intersection, which are either strict transforms of
exceptional divisors, or strict transforms of curves of $\gp^2$
satisfying the above conditions 1) and 2) (the second condition is
obvious and the first one is a consequence of the adjunction
formula).




The sequence of dual cones $V_0^\vee \supset V_1^\vee \supset
\cdots$ is strictly decreasing and each cone $V_i^\vee$ contains
$P(Z_{\cf})$. Therefore, after repeating the steps 2, 3 and 4
finitely many times,  Condition (b) of Step 5 will not be
satisfied and, hence, the process described in the algorithm
stops.

Finally, our algorithm is justified by bearing in mind Proposition
\ref{nueva} and the fact that the final set $G$ is a maximal
subset of $\{Q_1,Q_2,\ldots\}$ among those whose elements are
invariant curves and $\{[\tilde{R}]\mid R\in G \}\cup
\{[\tilde{E}_q]\}_{q\in {\cal N}_{\cf}}$ is a linearly independent
subset of $A(Z_{\cf})$. Observe that all the curves in $G$ are
integral since they have been computed with degrees in increasing
order.\findemo \\

\begin{rem}
{\rm Polyhedrality of the cone $NE(Z_{\cf})$ ensures that
Algorithm \ref{alg3} ends. Notwithstanding, one can apply it
anyway and, if it stops after a finite number of steps, one gets
an output as it is stated in the  algorithm. When  an independent
system of algebraic solutions $S$ is obtained, one might decide
whether, or not, $\cf$ admits a rational first integral either by
using Algorithm \ref{alg2} (when some of the conditions (1), (2)
and (3) stated before it held) or by applying results from Section
\ref{DOS}, as Proposition \ref{memo} or Remark 1. Let us see an
example.

}
\end{rem}

\begin{exa}
{\rm

Let $\cf$ be the foliation on $\gp^2_{\gc}$ defined by the
projective 1-form $\Omega=AdX+BdY+CdZ$, where:
$$A=-3X^2Y^3+9X^2Y^2Z-9X^2YZ^2+3X^2Z^3, \;\; B=3X^3Y^2-6X^3YZ-5Y^4Z+3X^3Z^2
\;\;\mbox{and}$$
$$C=-3X^3Y^2+5Y^5+6X^3YZ-3X^3Z^2.$$
The configuration $\ck_{\cf}$ consists of the union of two chains
$\{q_i\}_{i=1}^{19}\cup \{q_i\}_{i=20}^{23}$ with the following
additional proximity relations: $q_4$ is proximate to $q_2$,
$q_{22}$  to $q_{20}$ and $q_{23}$ to $q_{21}$. Since the unique
dicritical exceptional divisor is $E_{q_{19}}$, we have that ${\cal
B}_{\cf}=\{q_i\}_{i=1}^{19}$ and ${\cal
N}_{\cf}=\{q_i\}_{i=1}^{18}$.

A priori, we do no know whether $NE(Z_{\cf})$ is, or not,
polyhedral. However, if we run Algorithm \ref{alg3}, it ends;
providing the independent system of algebraic solutions $S=\{C\}$,
where $C$ is the line defined by the equation $Y-Z=0$. Therefore,
$T_{\cf,S}=5L^*-2E_{q_1}^*-2E_{q_2}^*-\sum_{i=3}^{19} E_{q_i}^*$
and, since
$[T_{\cf,S}]=5[\tilde{C}]+3[\tilde{E}_{q_1}]+6[\tilde{E}_{q_2}]+10[\tilde{E}_{q_3}]+
\sum_{i=4}^{18} (19-i)[\tilde{E}_{q_i}]$,  Condition (2) before
Algorithm \ref{alg2} is satisfied and, hence, we can apply this
algorithm. The value $\Delta_{\cf}$ of Theorem \ref{gordo2} equals 1
and the algorithm returns the rational first integral $F_1/F_2$,
given by $F_1=Y^5-X^3Y^2+2X^3YZ-X^3Z^2$ and $F_2=(Y-Z)^5$. }
\end{exa}

The next proposition shows that {\it if the cone $NE(Z_{\cf})$ of
a foliation $\cf$ is polyhedral, then calling Algorithms
\ref{alg3} and \ref{alg2}, one can decide whether $\cf$ has a
rational first integral and, in the affirmative case, to compute
it.} The unique data we need in that procedure are the following:
a projective 1-form $\Omega$ defining $\cf$, the configuration of
dicritical points ${\cal B}_{\cf}$ and the set ${\cal N}_{\cf}$.

\begin{pro}

Let $\cf$ be a foliation on $\gp^2$ such that $NE(Z_{\cf})$ is a
polyhedral cone. Let $S$ be an independent system of algebraic
solutions obtained by calling Algorithm \ref{alg3}. Then, $S$
satisfies one of the conditions $(1)$, $(2)$ or $(3)$ described
before Algorithm \ref{alg2}.

\end{pro}

\noindent {\it Proof}. Assume that Condition $(1)$ does not hold.
Let $V$ be the convex cone obtained by calling Algorithm \ref{alg3}.
Since $y^2\geq 0$ for all $y\in V^\vee$, a similar reasoning to that
given in the proof of Proposition \ref{nueva} (taking $T_{\cf,S}$
instead of $D_{\cf}$ and $V$ in place of $\con(W)$) shows that
$V\cap [T_{\cf,S}]^\bot$ is a face of $V$ which contains
$[T_{\cf,S}]$ and that, for all $x\in [T_{\cf,S}]^\bot\setminus
\{0\}$, $x^2\leq 0$ and so $x^2<0$ if, and only if, $x$ is not a
(real) multiple of $[T_{\cf,S}]$.

From the above facts, it is straightforward to deduce that the class
$[T_{\cf,S}]$ does not belong to any proper face of $V\cap
[T_{\cf,S}]^\bot$. Thus  $[T_{\cf,S}]$ is in  the relative interior
of $V\cap [T_{\cf,S}]^\bot$, since each non-zero element of a
polyhedral convex cone belongs to the relative interior of a unique
face. Then, if $V\cap [T_{\cf,S}]^\bot$ is a simplicial cone, it is
clear that Condition $(2)$ given before Algorithm \ref{alg2} is
satisfied. Otherwise, $[T_{\cf,S}]$ admits, at least, two different
decompositions as linear combination (with rational positive
coefficients) of classes of irreducible curves on $Z_{\cf}$
belonging to $V\cap [T_{\cf,S}]^\bot$. Therefore,
$h^0(\gp^2,{\pi_{\cf}}_*\co_{Z_{\cf}}(\lambda T_{{\cal F},S}))\geq
2$ for some positive integer $\lambda$ and Condition $(3)$ holds.\findemo\\


 In \cite{gal1},  we gave conditions which imply the
polyhedrality of the cone of curves of the smooth projective
rational surface obtained by blowing-up at the infinitely near
points of a configuration $\cal{C}$ over $\gp^2$. One of them only
depends  on the proximity relations among the points of $\cal{C}$
and it holds for a wide range of surfaces whose anticanonical
bundle is not ample. Now, we shall recall this condition, but
first we introduce some necessary  notations.

Consider the morphism $f:Z\rightarrow \gp^2$ given by the
composition of the blow-ups centered at the points of ${\cal C}$ (in
a suitable order). For every $p\in {\cal C}$, the exceptional
divisor $E_p$ obtained in the blow-up centered at $p$ defines a
valuation of the fraction field of $\co_{\gp^2, p_0}$, where $p_0$
is the image of $p$ on $\gp^2$ by the composition of blowing-ups
which allow to obtain $p$. So, $p$  is associated with a simple
complete ideal $I_p$ of the ring $\co_{\gp^2, p_0}$ \cite{lip}.
Denote by $D(p)$ the exceptionally supported divisor of $Z$ such
that $I_p \co_Z=\co_Z(-D(p))$ and by $K_Z$ a canonical divisor on
$Z$. Let $G_{\cal C}=(g_{p,q})_{p,q\in {\cal C}}$ be the square
symmetric matrix defined as follows: $g_{p,q}=-9 D(p)\cdot
D(q)-\bigl(K_Z \cdot D(p)\bigr)\bigl(K_Z\cdot D(q)\bigr)$ (see
\cite{gal1}, for an explicit description of $G_{\cal C}$ in terms of
the proximity graph of ${\cal C}$).

\begin{de}
\label{suficiente} {\rm A configuration ${\cal C}$ as above is
called to be {\it P-sufficient} if  ${\bf x} G_{\cal C} {\bf
x}^t>0$ for all non-zero vector ${\bf x}=(x_1,x_2,\ldots,x_r)\in
\gr^r$ with non-negative coordinates (where $r$ denotes the
cardinality of ${\cal C}$).}

\end{de}

This definition gives the cited condition, thus  {\it if ${\cal
C}$ is a P-sufficient configuration, then the cone of curves of
the surface $Z$ obtained by blowing-up their points is polyhedral}
\cite[Th. 2]{gal1}. Notice that using the criterion given in
\cite{gad}, it is possible to decide whether a configuration is
P-sufficient. Moreover, when the configuration ${\cal C}$ is a
{\it chain} (that is, the partial ordering $\geq $ defined in
Section \ref{UNO} is a total ordering), a very simple to verify
criterion can be given: {\it ${\cal C}$ is P-sufficient if the
last entry of the matrix $G_{\cal C}$ is positive} \cite[Prop.
6]{gal1}. Also, it is worthwhile to add that whichever
configuration of cardinality less than 9 is P-sufficient.

Taking into account the above facts and applying \cite[Lem.
3]{gal3}, it can be proved the following result which shows that,
if the configuration of dicritical points of an arbitrary
foliation $\cf$ on $\gp^2$ is P-sufficient, then Algorithm
\ref{alg3} and Proposition \ref{memo} can be applied to decide
whether $\cf$ has a rational first integral (and to compute it).

\begin{pro}
\label{lacinco}
 Let ${\cal F}$ be a foliation on $\gp^2$ such that
the configuration ${\cal B}_{\cf}$ is P-sufficient. Then:

\begin{itemize}
\item[(a)] The cone $NE(Z_{\cf})$ is polyhedral.

\item[(b)]If $S$ is an independent system of algebraic solutions
obtained by calling Algorithm \ref{alg3}, then either
$T_{\cf,S}^2\not=0$,  $T_{{\cal F},S}\cdot \tilde{E}_q<0$ for some
$q\in {\cal B}_{\cal F}\setminus {\cal N}_{\cal F}$ or
$K_{Z_{\cf}}\cdot T_{\cf,S}<0$.
\end{itemize}

\end{pro}



To end this paper, we give two illustrative examples where we have
applied the above results and algorithms.

\begin{exa}
{\rm Let $\{\cf_a\}_{a\in \gc}$ be  the one-parameter family of
foliations on the projective plane over the complex numbers $\gc$
defined by the projective 1-form $\Omega=AdX+BdY+CdZ$, where:
$$A=Z(aXZ-Y^2+Z^2), \;\;\;
B=Z(X^2-Z^2)\;\;\; \mbox{and} \;\;\; C=XY^2-aX^2Z-XZ^2-X^2Y+YZ^2.$$
Set $Q:=\{\frac{q^2-p^2}{q^2}|\; p,q \in \mathbb{Z},\; q \neq 0\}$
and consider the following points on $\gp^2_\gc$: $U=(1:0:0)$,
$W=(0:1:0)$, $R_a= \left(1:\sqrt{1+a}:1\right)$ when $-a \in Q$ and
$S_a=\left(1:\sqrt{1-a}:-1\right)$ whenever $a \in Q$. In order to
determine the foliations of the family $\{\cf_a\}_{a\in \gc}$ with a
(or without any) rational first integral, we shall
distinguish three cases:\\

\noindent {\it Case 1}: $a$ and $-a$ are not in $Q$.\\

\noindent In this case, ${\cal K}_{\cf_a}=\{U,W\}$, ${\cal
B}_{\cf_a}=\{W\}$ and ${\cal N}_{\cf_a}=\emptyset$.
Since the cardinality of ${\cal B}_{\cf_a}$ is 1, the inequality
$d(X/Z)\wedge \Omega \neq 0$ shows that $\cf_a$ has no rational
first integral.\\

\noindent {\it Case 2:}  $a \neq 0$ and either $a$ or $-a$ belong to $Q$.\\

\noindent Now ${\cal K}_{\cf_a}=\{U,W \} \cup \mathfrak{D}_a$ and
${\cal B}_{{\cal F}_a}=\{W\}\cup \mathfrak{D}_a$, where
$\mathfrak{D}_a$ stands for a configuration whose unique point on
$\gp^2$ is $S_a$ (respectively, $R_a$) whenever $a \in Q$
(respectively, $-a \in Q$). $\mathfrak{D}_a$ strongly depends on the
value $a$; we shall illustrate it taking two specific values for $a$
and running our algorithms in order to decide about the existence of
a rational first integral in each case.

If $a=5/9$, $\mathfrak{D}_a$ consists of a chain of 3 points
$\{q_1,q_2,q_3\}$, where $q_3$ is proximate to $q_1$ and provides
the unique dicritical exceptional divisor associated with points in
$\mathfrak{D}_a$. Since the configuration ${\cal B}_{\cf_a}$ is
P-sufficient, the cone $NE(Z_{\cal F})$ is polyhedral and if we run
Algorithm \ref{alg3} for the foliation $\cf_a$, it will end. The
sequence of convex cones $V_0\subset V_1$ computed by the algorithm
(see the explanation in page \pageref{explanation}) is the
following: $V_0=\con(\{[\tilde{E}_{W}], [\tilde{E}_{q_1}],
[\tilde{E}_{q_2}], [\tilde{E}_{q_3}]\})$ and $V_1=\con(V_{0}\cup
\{[L^*-{E}_{W}^*-{E}_{q_1}^*-{E}_{q_2}^*]\})$. The dual cone of
$V_1$ is spanned by the classes $[L^*]$, $[L^*-E_{W}^*]$,
$[L^*-E_{q_1}^*]$, $[2L^*-E_{q_1}^*-E_{q_2}^*]$ and
$[3L^*-2E_{q_1}^*-E_{q_2}^*-E_{q_3}^*]$. Since all these generators
have non-negative self-intersection, (b) of Step 5 of the above
mentioned algorithm is not satisfied for $V_1$ and, then, the
algorithm ends. The unique element of the obtained set $G$ is the
line with equation $X+Z=0$, so Step 6 returns ``0'' and, therefore,
${\cal F}_a$ has no rational first integral.

If $a=- \frac{861}{100}$, $\mathfrak{D}_a$ is a chain of 13 points
$\{q_i\}_{i=1}^{13}$, characterized by the fact that $q_{13}$ is
proximate to $q_{3}$, and ${\cal N}_{\cf}=\{q_i\}_{i=1}^{12}$. Since
the configuration ${\cal B}_{\cf_a}$ is also P-sufficient, Algorithm
\ref{alg3} will end. The sequence of obtained convex cones is
$V_0\subset V_1\subset V_2$, where
$V_0=\con(\{[\tilde{E}_{q}]\}_{q\in {\cal B}_{\cf_a}})$,
$V_1=\con(V_0\cup \{ [L^*-E_W^*-E_{q_1}^*]\})$ and $V_2=\con(V_1\cup
\{ [L^*-E_{q_1}^*-E_{q_2}^*]\})$. The dual cone of $V_2$ has 27
extremal rays which are spanned by classes with non-negative
self-intersection. Moreover, the  set $G$ is the same as above. So,
Algorithm \ref{alg3} returns ``0'' and, thus, ${\cal F}_a$ has no
rational first integral.\\

It is possible to discard the existence of a rational first integral
for any foliation $\cf_a$ corresponding to the current case. Indeed,
reasoning by contradiction, assume that $\cf_a$ has a rational first
integral and consider the algebraic solution $C$ of $\cf_a$ defined
by the line with equation $Z=0$. Since the self-intersection of the
class of its strict transform on $Z_{\cf_a}$ vanishes, by Theorem
\ref{teorema} it is easy to deduce that $D_{\cf_a}=L^*-E_W^*$. The
set $\{X,Z\}$ is a basis of
$H^0(\gp^2,{\pi_{\cf_a}}_*\co_{Z_{\cf_a}}(D_{\cf_a}))$ and, since
$d(X/Z)\wedge \Omega\not=0$, we get a contradiction.\\

\noindent {\it Case 3:} $a=0$\\

\noindent Here, $ {\ck}_{\cf_0} = {\cal B}_{\cf_0}=\{R=R_0,
S=S_0,U,W\}$, ${\cal N}_{\cf_0}=\emptyset$. Applying Algorithm
\ref{alg3} we get the following independent system of algebraic
solutions for $\cf_0$:
$$G=\{l_{R,S}, l_{S,U}, l_{U,W}, l_{S,W}\},$$ where $l_{p,q}$
denotes the line joining $p$ and $q$ ($p, q \in {\cal B}_{\cf_0}$).
Its associated divisor $T_{\cf_0,G}$ is
$2L^*-E^*_R-E^*_S-E^*_U-E^*_W$. By Proposition \ref{memo}, we get
that if $\cf_0$ admits rational first integral, then $D_{\cf_0}$
should coincide with $T_{\cf_0,G}$. The space of global sections
$H^0(\gp^2,{\pi_{\cf}}_*\co_{Z_{\cf}}(T_{\cf_0,G}))$ has dimension 2
and it is generated by $F_1=(X+Z)(Z-Y)$ and $F_2=Z(Y-X)$. So we
conclude, after checking it, that the rational function $F_1/F_2$ is
a rational first integral of ${\cf}_0$. }
\end{exa}

\begin{figure}[b]\label{fig5}
\setlength{\unitlength}{1mm}
\begin{center}
\begin{picture}(45,25)

\put(0,0){\circle*{2}} \put(0,0){\line(0,1){10}} \put(-5,0){$q_1$}

\put(0,10){\circle*{2}} \put(-5,10){$q_2$}

\put(40,0){\circle*{2}} \put(40,0){\line(0,1){10}}
\put(45,0){$q_7$}

\put(40,10){\circle*{2}}\put(40,10){\line(0,1){10}}
\put(45,10){$q_8$}

\put(40,20){\circle*{2}} \put(45,20){$q_9$}

\qbezier[50](40,0)(35,10)(40,20)

\put(20,0){\circle*{2}} \put(25,0){$q_3$}

\put(10,10){\circle*{2}} \put(20,0){\line(0,1){10}}
\put(5,10){$q_4$}

\put(20,10){\circle*{2}} \put(20,0){\line(-1,1){10}}
\put(15,10){$q_5$}

\put(30,10){\circle*{2}} \put(20,0){\line(1,1){10}}
\put(25,10){$q_6$}

\end{picture}
\end{center}
\caption{The proximity graph of ${\cal B}_{\cf}$ in Example
\ref{ex5}}
\end{figure}

\begin{exa}\label{ex5}
{\rm

Now, let ${\cal F}$ be a foliation as above defined by the 1-form
$\Omega=AdX+BdY+CdZ$, where
$$A=X^4Y^3Z+5X^3Y^4Z+9X^2Y^5Z+7XY^6Z+2Y^7Z+X^4Z^4-X^3YZ^4,$$
$$B=-3X^5Y^2Z-13X^4Y^3Z-21X^3Y^4Z-15X^2Y^5Z-4XY^6Z+2X^4Z^4 \;\; \mbox{and}$$
$$C=2X^5Y^3+8X^4Y^4+12X^3Y^5+8X^2Y^6+2XY^7-X^5Z^3-X^4YZ^3.$$
Resolving $\cf$, we get ${\ck}_{\cf}={\cal B}_{\cf}=
\{q_i\}_{i=1}^{9}$. Its proximity graph is given in Figure 3. ${\cal
B}_{\cf}$ is P-sufficient (use the above mentioned result given in
\cite{gad}) and ${\cal N}_{\cf}= \{q_1,q_3,q_7,q_8\}$. The output
$G$ of Algorithm \ref{alg3}  is given by the curves with equations
$X=0$, $X+Y=0$, $Z=0$, $XY+Y^2+XZ=0$ and $jXY+jY^2+XZ$, where $j$ is
a primitive cubic  root of unity. So, $G$ is an independent system
of algebraic solutions. Its associated divisor $T_{\cf,G}$ is
$$
6L^*-3 \sum_{i=1}^3 E^*_{q_i}- \sum_{i=4}^6 E^*_{q_i}- 2E^*_{q_7}-
\sum_{i=8}^9 E^*_{q_i}.
$$
By Proposition \ref{memo}, if ${\cal F}$ has a rational first
integral, then the pencil $\cp_{\cf}$ is $\gp
H^0(\gp^2,{\pi_{\cf}}_*\co_{Z_{\cf}}(T_{\cf,G}))$. A basis of this
projective space is given by $F_1=(X+Y)^2 X^2 Z^2$ and $F_2=(X+Y)^3
Y^3 + X^3 Z^3$. Finally, the equality $d(F_1/F_2)\wedge \Omega =0$
shows that $F_1/F_2$ is a rational first integral of $\cf$.

}
\end{exa}

\vspace{3mm}
\par
\begin{center}
\footnotesize Dept. de Matemàtiques (ESTCE), UJI, Campus Riu Sec. \\
\footnotesize 12071  Castelló. SPAIN. \\ \footnotesize
galindo@mat.uji.es; monserra@mat.uji.es
\end{center}


\begin{thebibliography}{99}
\footnotesize \setlength{\baselineskip}{3mm}
\bibitem{aut} L. Autonne, Sur la th\'eorie des \'equations
diff\'erentielles du premier ordre et du premier degr\'e, {\em J.
\'Ecole Polytech.} {\bf 61} (1891), 35---122; {\bf 62} (1892),
47---180. \vspace{-1.5mm} \bibitem{bau} N.N. Bautin, On periodic
solutions of a system of differential equations, {\it Akad. Nauk.
SSSR. Prikl. Mat. Meh.} {\bf 18} (1954) 128, (in Russian).
\vspace{-1.5mm}
\bibitem{beau} A. Beauville, {\it Complex algebraic surfaces}, London Math. Soc. Student Texts {\bf 34}, Cambridge
University Press, 1996. \vspace{-1.5mm}
\bibitem{be} I. Bendixson, Sur les points singuliers d'une
\'equation diff\'erentielle lin\'eaire, {\it Ofv. Kongl. Ventenskaps
Akademiens Forhandlinger} {\bf 148} (1895), 81---89. \vspace{-1.5mm}
\bibitem{brun} M. Brunella, {\it Birational Geometry of
Foliations}, Springer, 2000. \vspace{-1.5mm}
\bibitem{ca-ca} A. Campillo and M. Carnicer,
Proximity inequalities and bounds for the degree of invariant curves
by foliations of $\gp_{\mathbb{C}}^2$, {\em Trans. Amer. Math. Soc.}
{\bf 349 (9)} (1997), 2211---2228. \vspace{-1.5mm}
\bibitem{cam}  A. Campillo, G. González-Sprinberg and M.
Lejeune-Jalabert, Clusters of infinitely near points, {\it Math.
Ann.} {\bf 306 (1)}  (1996), 169---194.
\vspace{-1.5mm}
\bibitem{cam-gon}  A. Campillo and G. González-Sprinberg, On characteristic
cones, clusters and chains of infinitely near points, in Brieskorn
Conference Volume, Progr. Math {\bf 168}, Birh\"auser, 1998.
\vspace{-1.5mm}
\bibitem{camp} A. Campillo, O. Piltant and  A. Reguera,
Cones of curves and of line bundles on surfaces associated with
curves having one place at infinity, {\em Proc. London Math. Soc.}
{\bf 84} (2002), 559---580.\vspace{-1.5mm}
\bibitem{car} M. Carnicer, The Poincar\'{e}
problem in the nondicritical case, {\em Ann. Math.} {\bf 140}
(1994), 289---294. \vspace{-1.5mm}
\bibitem{casas} E. Casas-Alvero, {\it Singularities of plane
curves}, London Math. Soc. Lecture Note Series {\bf 276}, Cambridge
University Press, 2000. \vspace{-1.5mm}
\bibitem{ch-gi} J. Chavarriga, H. Giacomini and J. Gin\'e, An
improvement to Darboux integrability theorem for systems having a
center, {\it Appl. Math. Lett.} {\bf 12} (1999), 85---89.
\vspace{-1.5mm}
\bibitem{cha-lli} J. Chavarriga and J. Llibre, Invariant
algebraic curves and rational first integrals for planar
polynomial vector fields, {\it J. Diff. Eq.} {\bf 169 (1)} (2001),
1---16. \vspace{-1.5mm}
\bibitem{cut} S.D. Cutkosky and V. Srinivas, On a problem of  Zariski
on dimension of linear systems, {\it Ann. Math.} {\bf 137} (1993),
531---559. \vspace{-1.5mm}
\bibitem{dar} G. Darboux, M\'emoire sur les \'equations
diff\'erentielles alg\'ebriques du premier ordre et du premier
degr\'e (M\'elanges), {\em Bull. Sci. Math.} {\bf 32} (1878),
60---96; 123---144; 151---200. \vspace{-1.5mm}
\bibitem{do-lo} V.A.
Dobrovol'skii, N.V. Lokot and J.M. Strelcyn, Mikhail Nikolaevich
Lagutinskii (1871-1915): Un math\'ematicien m\'econnu, {\it Historia
Mathematica} {\bf 25} (1998), 245---264. \vspace{-1.5mm}
\bibitem{eis} D. Eisenbud, M. Green and J. Harris, Cayley-Bacharach
theorems and conjectures, {\it Bull. Amer. Math. Soc. (N.S.)} {\bf
33 (3)} (1996), 295---324.\vspace{-1.5mm}
\bibitem{es-kl} E. Esteves  and S. Kleiman, Bounds on leaves of one-dimensional foliations.
 Preprint 2002, math AG 0209113.
\vspace{-1.5mm}
\bibitem{gad} J. W. Gaddum, Linear inequalities and quadratic
forms, {\it Pacific J. Math.} {\bf 8} (1958), 411---414.
\vspace{-1.5mm}
\bibitem{gal1} C. Galindo and  F. Monserrat, The cone of curves
associated to a plane configuration, {\it Comm. Math. Helv.} {\bf
80} (2005), 75---93. \vspace{-1.5mm}
\bibitem{gal3} C. Galindo and F. Monserrat, The total coordinate ring of a smooth projective surface,
 {\it J. Algebra} {\bf 284} (2005), 91---101.
 \vspace{-3.5mm}
\bibitem{julio} J. Garcia-de la Fuente, Geometr\'{\i}a de los
sistemas lineales de series de potencias en dos variables, Ph. D.
thesis, Valladolid University (1989), (in Spanish). \vspace{-1.5mm}
\bibitem{G-M} X. G\'{o}mez-Mont and  L. Ortiz, {\it Sistemas din\'amicos holomorfos en
superficies}, Aportaciones Matem\'aticas {\bf 3}, Sociedad
Matemática Mexicana, 1989 (in Spanish). \vspace{-1.5mm}
\bibitem{hart} R. Hartshorne, {\it Algebraic Geometry},
Springer-Verlag, 1977. \vspace{-1.5mm} \bibitem{hew} C. Hewitt,
Algebraic invariant curves in cosmological dynamical systems and
exact solutions. {\it Gen. Relativity Gravitation} {\bf 23}
(1991), 1363---1384. \vspace{-1.5mm} \bibitem{hil} D. Hilbert,
Mathematische problem (lecture), Second Internat. Congress Math.
Paris, Nachr. Ges. Wiss. G\"ottingen Math-Phys. Kl. 1900, 253-297;
English translation {\it Bull. Amer. Math. Soc.} {\bf 8} (1902),
437-479; reprinted in Mathematical developments arising from
Hilbert problems, Proc. Sympos. Pure Math. {\bf 28} (1976),
1---34. \vspace{-1.5mm}
\bibitem{jou} J.P. Jouanolou, Equations de Pfaff
algébriques, Lecture Notes in Mathematics {\bf 708},
Springer-Verlag, 1979. \vspace{-1.5mm}
\bibitem{kaliman} S. Kaliman, Two remarks on polynomials in two
variables, {\it Pacific J. Math.} {\bf 154} (1992), no. 2,
285---295.\vspace{-1.5mm}
\bibitem{kle} S. Kleiman, Towards a numerical theory of ampleness,
{\it Ann. Math.} {\bf 84} (1966), 293---349. \vspace{-1.5mm}
\bibitem{kollar} J. Koll\'ar and S. Mori,
{\it Birational geometry of rational varieties}, Cambridge Tracts in
Mathematics {\bf 134}, Cambridge University Press, 1998.
\vspace{-1.5mm}
\bibitem{le} S. Lefschetz, On a theorem of Bendixson, {\it J.
Diff. Eq.} {\bf 4} (1968), 66---101. \vspace{-1.5mm}
\bibitem{l-n} A. Lins-Neto,
Some examples for the Poincar\'e and Painlev\'e problems, {\it Ann.
Sc. \'Ec. Norm. Sup.} {\bf 35} (2002), 231---266. \vspace{-3.5mm}
\bibitem{lip} J. Lipman, Rational singularities with applications to
algebraic surfaces and unique factorization. {\it Publ. Math. IHES}
{\bf 36} (1969), 195---279.\vspace{-1.5mm}
\bibitem{lli}  J. Llibre and G. \'Swirszcz,  Relationships
between limit cycles and algebraic invariant curves for quadratic
systems. Preprint 2003. \vspace{-1.5mm}
\bibitem{man1} Y. Man, Computing closed form solutions of first
order ODEs with the Prelle-Singer procedure, {\it J. Symb. Comp.}
{\bf 16} (1993), 423---443. \vspace{-1.5mm}
\bibitem{man2} Y. Man and A. MacCallum, A rational approach to the
Prelle-Singer algorithm, {\it J. Symb. Comp.} {\bf 24} (1997),
31---43. \vspace{-1.5mm}
\bibitem{man} Y. Manin, {\it Cubic forms. Algebra, Geometry,
Arithmetic}, North Holland, 1974. \vspace{-1.5mm}
 \bibitem{mor} S. Mori, Threefolds whose
canonical bundles are not numerically effective, {\it Ann. Math.}
{\bf 116} (1982), 133---176. \vspace{-1.5mm}
\bibitem{oda} T. Oda, {\it Convex bodies and algebraic geometry,
an introduction to the theory of toric varieties}, Ergebnisse der
Math. {\bf 15}, Springer-Verlag, 1988.\vspace{-1.5mm}
\bibitem{pai} P. Painlev\'e,
 ``Sur les int\'egrales alg\'ebriques des
\'equations diff\'erentielles du premier ordre" and ``M\'emoire sur
les \'equations diff\'erentielles du premier ordre" in Ouvres de
Paul Painlev\'e, Tome II, \'Editions du Centre National de la
Recherche Scientifique 15, quai Anatole-France, Paris 1974.
\vspace{-1.5mm} \bibitem{per} J.V. Pereira, Vector fields, invariant
varieties and linear systems, {\it Ann. Inst. Fourier} {\bf 51(5)}
(2001), 1385---1405. \vspace{-1.5mm}
\bibitem{poi1} H. Poincar\'e, M\'emoire sur les courbes
d\'efinies par les \'equations diff\'erentielles, {\em J. Math.
Pures Appl.}  {\bf 3 (7)} (1881), 375---442; {\bf 3 (8)} (1882),
251---296; {\bf 4 (1)}  (1885), 167---244; in Oeuvres de Henri
Poincar\'e, vol. I, Gauthier-Villars, Paris 1951, 3---84, 95---114.
\vspace{-1.5mm}
\bibitem{poi2} H. Poincar\'e, Sur l'int\'egration alg\'ebrique des
 \'equations diff\'erentielles du premier ordre et
du premier degr\'e (I), {\em Rend.  Circ. Mat. Palermo} {\bf 5}
(1891), 161--191. \vspace{-1.5mm}
\bibitem{poi3} H. Poincar\'e, Sur l'int\'egration alg\'ebrique des
\'equations diff\'erentielles du premier ordre et du premier degr\'e
(II), {\em Rend.  Circ. Mat. Palermo} {\bf 11} (1897), 193---239.
\vspace{-1.5mm}
\bibitem{pr-si} M.J. Prelle and
M.F. Singer, Elementary first integrals of differential equations,
{\it Trans. Amer. Math. Soc.} {\bf 279} (1983), 215---229.
\vspace{-1.5mm}
\bibitem{rock} R.T. Rockafellar, {\it Convex Analysis}, Princeton
Univ. Press, 1970. \vspace{-1.5mm} \bibitem{sch} D. Scholomiuk,
Algebraic particular integrals, integrability and the problem of the
centre, {\it Trans. Amer. Math. Soc.} {\bf 338} (1993), 799---841.
\vspace{-1.5mm}
\bibitem{seid}  A. Seidenberg, Reduction of singularities of
the differentiable equation $Ady=Bdx$, {\it Amer. J. Math.} {\bf 90}
(1968), 248---269. \vspace{-1.5mm}
\bibitem{sma} S. Smale, Mathematical problems
for the next century, {\em Math. Intelligencer} {\bf 20} (1998),
7---15. \vspace{-1.5mm}
\bibitem{zar} O. Zariski, The theorem of Riemann-Roch for high multiples
of an effective divisor on an algebraic surface, {\it Ann. Math.}
{\bf 76 (3)} (1962), 560---615.
\end{thebibliography}
\end{document}